\documentclass[a4paper,11pt]{article}
\usepackage[top=25mm,left=22mm]{geometry}
\usepackage{amsmath,amssymb,amsthm,amsfonts}
\usepackage{lineno,url}
\usepackage{graphicx}
\usepackage{multirow}
\usepackage{pstricks}
\usepackage{paralist, longtable}
\newtheorem{theorem}{Theorem}[section]
\newtheorem{remark}[theorem]{Remark}
\def\bexe{\begin{exercise}}\def\eexe{\eex\end{exercise}}
\def\bsol{\begin{solution}}\def\esol{\eex\end{solution}}
\def\bexa{\begin{example}}\def\eexa{\end{example}}
\def\brem{\begin{remark}}\def\erem{\end{remark}}
\def\bthm{\begin{theorem}}\def\ethm{\end{theorem}}
\def\blem{\begin{lemma}}\def\elem{\end{lemma}}
\def\bcor{\begin{corollary}}\def\ecor{\end{corollary}}
\def\bdefi{\begin{definition}}\def\edefi{\end{definition}}
\newcommand{\IDEA}{\textbf{Idea of the Proof.} }

\def\bmip{\begin{minipage}{\textwidth}}\def\emip{\end{minipage}}
\def\huga#1{\begin{gather} #1 \end{gather}}
\def\hugast#1{\begin{gather*} #1 \end{gather*}}
\def\hual#1{\begin{align} #1 \end{align}}

\newcommand{\R}{{\mathbb R}}
\newcommand{\N}{{\mathbb N}}
\newcommand{\Z}{{\mathbb Z}}

\def\CG{{\cal G}}
\def\CO{{\cal O}}

\def\per{{\rm per}}

\def\ga{\gamma}\def\om{\omega}
\def\noi{\noindent}
\def\pa{{\partial}}\def\lam{\lambda}
\newcommand{\bi}{\begin{itemize}}\newcommand{\ei}{\end{itemize}}
\newcommand{\ben}{\begin{enumerate}}\newcommand{\een}{\end{enumerate}}
\newcommand{\bce}{\begin{center}}\newcommand{\ece}{\end{center}}
\newcommand{\bci}{\begin{compactitem}}\newcommand{\eci}{\end{compactitem}}
\newcommand{\bcen}{\begin{compactenum}}\newcommand{\ecen}{\end{compactenum}}
\newcommand{\reff}[1]{(\ref{#1})}

\newcommand{\spr}[1]{\left\langle #1 \right\rangle}
\newcommand{\vs}[1]{{\vspace{#1}}}

\def\eps{\varepsilon}

\newcommand{\barr}{\begin{array}}\newcommand{\earr}{\end{array}}
\newcommand{\bpm}{\begin{pmatrix}}\newcommand{\epm}{\end{pmatrix}}
\newcommand{\bsm}{\left(\begin{smallmatrix}}
\newcommand{\esm}{\end{smallmatrix}\right)}
\newcommand{\ba}{\begin{array}}\newcommand{\ea}{\end{array}}
\def\dd{\, {\rm d}}\def\ri{{\rm i}}

\def\er{{\rm e}}

\def\om{\omega}\def\Om{\Omega}

\def\del{\delta}

\def\eex{\hfill\mbox{$\rfloor$}}

\def\Del{\Delta}

\def\sig{\sigma}
\def\al{\alpha}

\def\Lam{\Lambda}

\def\bd{\begin{displaymath}} \def\ed{\end{displaymath}}
\def\ba{\begin{array}} \def\ea{\end{array}}
  
\def\eps{\varepsilon}

\def\ig{\includegraphics}\def\pdep{{\tt pde2path}}
\def\pdepb{{\tt p2p2}}
\def\mlab{{\tt Matlab}}\def\ptool{{\tt pdetoolbox}}

\def\bcs{BC}\def\per{{\rm per}}\def\tw{\tilde w}

\def\ssmatrix#1{\bigl(\begin{smallmatrix}#1 \end{smallmatrix} \bigr)}
\def\wact{w_{{\rm act}}}\def\wt{\tilde{w}}
\renewcommand{\arraystretch}{1}\renewcommand{\baselinestretch}{1.0}
\def\medskip{}\def\bigskip{}
\def\DT{\Delta T}\def\uP{u_1}\def\uM{u_2}

\begin{document}
\mbox{}\vspace{0.1cm}\begin{center}\Large
pde2path - version 2.0: faster FEM, 
multi-parameter continuation, nonlinear boundary conditions, 
and periodic domains -- a short manual\\[2mm] 
\normalsize
Tomas Dohnal$^1$, Jens D.M.\ Rademacher$^2$, Hannes Uecker$^3$, Daniel Wetzel$^4$,  \\[2mm]
\footnotesize
$^1$ Fakult\"at f\"ur Mathematik, TU Dortmund, D44227 Dortmund, 
dohnal@mathematik.tu-dortmund.de\\
$^2$ Fachbereich Mathematik, Universit\"at Bremen, D28359 Bremen, 
jdmr@uni-bremen.de \\
$^3$ Institut f\"ur Mathematik, Universit\"at Oldenburg, D26111 Oldenburg, 
hannes.uecker@uni-oldenburg.de\\
$^4$  Institut f\"ur Mathematik, Universit\"at Oldenburg, D26111 Oldenburg, 
daniel.wetzel@uni-oldenburg.de \\[2mm]
\normalsize
\today
\end{center}
\begin{abstract}\noindent
\pdep\ 2.0 is an upgrade of the continuation/bifurcation package \pdep\ for 
elliptic systems of PDEs over bounded 2D domains, based on \mlab's \ptool. 
The new features include a more efficient use of FEM, 
easier switching between different 
single parameter continuations, genuine multi--parameter 
continuation (e.g., fold continuation), more efficient implementation 
of nonlinear boundary conditions, cylinder and 
torus geometries (i.e., periodic boundary conditions), and a 
general interface for adding auxiliary equations like mass 
conservation or phase equations for continuation of traveling waves. 
The package (library, demos, manuals) can be downloaded at 
{\tt www.staff.uni-oldenburg.de/hannes.uecker/pde2path}. 
\end{abstract}
\tableofcontents 

\section{Introduction}\label{i-sec}
\pdep, based on the FEM of the \mlab\ \ptool, is a continuation/bifurcation 
package for elliptic systems of PDEs of the form 
\huga{\label{gform}
G(u,\lam):=-\nabla\cdot(c\otimes\nabla u)+a u-b\otimes\nabla u-f=0, 
}
where $u=u(x)\in\R^N$, $x\in\Omega\subset\R^2$ some 
bounded domain, $\lam\in\R^p$ is a parameter (vector), 
$c\in\R^{N\times N\times 2\times 2}$,
$b\in\R^{N\times N\times 2}$, 
 $a\in\R^{N\times N}$ and $f\in\R^N$
 can  depend on $x,u,\nabla u$, and, of course, 
parameters. The boundary conditions (\bcs) are 
``generalized Neumann'' of the form
\begin{align}\label{gnbc}
{\bf n}\cdot (c \otimes\nabla u) + q u = g,
\end{align}
where ${\bf n}$ is the outer normal and again $q\in \R^{N\times N}$
and $g\in \R^N$ may depend on $x$, $u$ and
parameters. These \bcs\ include prescribed flux \bcs, 
and a ``stiff spring'' approximation of Dirichlet \bcs\ via 
large prefactors in $q$ and $g$. 

For the basic ideas of continuation/bifurcation, 
the algorithms, and the 
class of systems we aim at, i.e., the meaning of the terms in \reff{gform} 
and the associated boundary conditions, 
we refer to \cite{p2p}, and the references therein. 
Here we explain a number of additional 
features in {\tt pde2path 2.0}, in short \pdepb, 
compared to the version documented in \cite{p2p}, 
and some changes in the underlying data structures. The new features 
include: 
\bcen
\item easy switching between different single parameter continuations; 
\item genuine multi--parameter continuation, in particular
automatic fold and branch point continuation;
\item general interface for adding auxiliary equations, such as mass 
conservation, or freezing-type equations for continuation 
of traveling waves; 
\item periodic domains: cylinder and torus geometries; 
\item fast FEM for a subclass of \reff{gform}, 
roughly where $c, a, b, q$ and $g$ are independent of $u$, 
i.e., where nonlinearity enters only through $f$; 
\item improved and more user-friendly plotting. 
\ecen

We explain these features by a number of examples, but
first describe the major structural changes (items 1,2,3). 
\brem\label{toolrem}{\rm The concept of \pdepb\ is that of a box of 
customizable tools. These tools (functions) are in {\tt p2p2/p2plib}, which 
must be in the \mlab\ path. When starting \mlab\ in the 
p2p2 home directory, execute {\tt setpde2path}. 
The demo directories are under 
{\tt p2p2/demos/}. Each demo (with name *) comes with one or more script-files 
{\tt *cmds.m}, which typically are organized in cells, i.e., 
should be stepped through cell by cell. To get help on any \pdepb\ function, 
e.g., {\tt cont}, type {\tt help cont} or {\tt doc cont}. 
To get started type {\tt help p2phelp}. Additional \mlab--internal 
and online html--help will be added shortly. 

To set up a new problem in \pdepb\ we recommend to copy a suitable 
demo directory (i.e., a demo directory which considers a similar problem) 
to a new working directory and start modifying the pertinent files. 
To customize any of the functions from {\tt p2p2/p2plib} we recommend 
to copy it to the working directory and modify it there 
(thus ``overloading'' the file 
from {\tt p2p2/p2lib}).}\eex 
\erem

\brem\label{comprem}{\rm The new data structure and different user interfaces 
mean that there is no downward compatibility with \cite{p2p}. 
On the other hand, we think that upgrading old \pdep\ files to \pdepb\ 
is quickly achieved, and that the data structure and user interfaces now 
have a final form.}\eex
\erem 

\paragraph{Parameters, auxiliary variables and  auxiliary equations.}
\label{ndsec}
A \pdepb\ problem is described by a matlab structure p.  The most drastic
change compared to \cite{p2p} is that no single distinguished
parameter $\lam$ appears in p anymore, but any number of auxiliary
variables, typically parameters, can be added. 
If the FEM mesh has $n_p=${\tt np} points and $N=${\tt neq} in
\reff{gform} we have $n_u=N n_p=${\tt p.nu=neq*p.np} unknown nodal values for $u$, 
(except in the case of periodic BC, see below), and
\texttt{p.u(1:p.nu)} contains these nodal values.  

The arbitrary number of auxiliary variables are stored in {\tt
  p.u(p.nu+1:end)} and can be ``passive'', serving as constant
parameters, or ``active'' unknowns to be solved for. In the following
we write, on the discrete level, $U=(u,w)=${\tt p.u}, where $u$
corresponds to (the nodal values of) the PDE variables in \reff{gform}
and $w$ the auxiliary variables.  Suppose there are $n_q+1$ active
variables $\wact\in\R^{n_q+1}$. Exactly one of these is the
``primary'' active parameter, and we write $\wact=(\wt,\al)$.  The
remaining $n_q$ active variables require $n_q$ additional
(`auxiliary') equations \huga{\label{qieq} q_i(U)=0, \quad
  i=1,\ldots,n_q. }  In the functions defining $G$ or its Jacobian a
typical first step is to split off the PDE part $u$ as shown in the
examples below.  The active auxiliary variables are selected by the
user in the array of indices {\tt p.nc.ilam}, whose first entry is the
primary continuation parameter. For different continuation tasks the
user may freely modify this list to choose different active, passive
and primary parameters. Thus, $\wact=(\wt,\al)$ is only a symbolic
notation, and the role of parameters (primary, active, passive) is
determined by {\tt p.nc.ilam}.  For convenience, the routines {\tt
  printaux} and {\tt getaux} can be used to obtain all or only the
active auxiliary variables. Internally, the routines {\tt au2u} and
{\tt u2au} are used to transform {\tt p.u} into the vector suitable
for the Newton-loop or back to the full {\tt p.u}.

Examples of additional equations are
\bci
\item prescribed mass: $\int u\dd x-m=0$, $m\in\R$, see, e.g., \S\ref{fchsec};
\item a phase condition $\langle \partial_x u, u-u_{\rm old}\rangle_2=0$ for the 
continuation of traveling waves, see, e.g., \S\ref{s:tw}. 
\eci

As discussed in \cite{p2p}, it is useful to give $u$ and the 
continuation parameter different weights in the arclength equation 
\huga{\label{e:arc}
p(U,s)=\spr{\dot U,U(s)-U_0}-(s-s_0)=0;
}
see \cite[\S2.1]{p2p}. In \pdepb~ this is extended to the active 
variables in  $\wact=(\wt,\al)$, and as scalar product in \reff{e:arc} we use, 
\huga{\label{xiqeq}
\spr{(u,\tw,\al),(v,\tilde{z},\beta)}:=\xi \spr{u,v}_2+\xi_q\spr{\tw,\tilde{z}}_2
+(1-(\xi+\xi_q)/2)\alpha\beta, 
}
with independent weights $\xi$ and $\xi_q$ and $\spr{\cdot,\cdot}_2$ 
the euclidean inner product. 

\medskip
\paragraph{Numerical approximation of $\partial_\lambda G$.}
In order to ease switching between different primary parameters, and since 
finite difference approximations of derivatives of $G$ with respect to just one 
parameter are relatively cheap, we deleted all explicit references to $\partial_\lam G$. Hence, the interfaces for the functions defining $G(u)$ and its Jacobian now read 
{\tt function [c,a,f,b]=G(p,u)} and {\tt  function [cj,aj,bj]=Gjac(p,u)},
see the examples below.

\paragraph{Substructures of {\tt p}: Names of numerical variables, switches, etc.}
In \pdepb\ the many switches and settings in the problem structure variable {\tt p} in 
\cite{p2p} are now grouped  as explained in Table \ref{tab1a}, i.e., 
function handles are entries in {\tt p.fuha}, numerical controls are 
entries in {\tt p.nc}, and so on. See also Appendix \ref{app1}. 
In particular, this makes it easier to get an overview over 
current parameter settings. For instance, to see the values of (all) the 
numerical control parameters for a given {\tt p}, type {\tt p.nc} on the 
command line. Of course, the user is free to add as many additional 
fields/variables to the structure {\tt p} as desired/needed. 
If there are many of these, then we recommend to organize them in a 
substructure {\tt p.usr}, say. 

An example of a (here predefined) 
``user''--field is {\tt p.usrlam}, which may contain target values for 
the primary parameter $\lam$. This means that if during 
continuation {\tt u(p.nu+p.nc.ilam(1))} passes a value $\lam^*$ 
in {\tt [p.usrlam,p.nc.lammin,p.nc.lammax]}, then the algorithm calculates and 
saves to file the solution at $\lam^*$. 
The names of files containing solution data at a continuation 
point have changed from the previous {\tt (p.dir/)p*.mat} to 
{\tt (p.dir/)pt*.mat}. Similarly, the solution at a bifurcation point is saved in {\tt bpt*.mat} and a fold point is saved in {\tt fpt*.mat}.

\LTcapwidth=\textwidth
\bce 
\begin{longtable}{| p{0.07\textwidth}|p{0.38\textwidth}|p{0.07\textwidth}|p{0.38\textwidth}|}
\caption{Fields in structure {\tt p}; see {\tt stanparam.m} in 
{\tt p2plib} for detailed information on the contents of these 
fields and the standard settings, and the reference card in Appendix 
\ref{app1}. \label{tab1a}}
\endfirsthead\endhead\endfoot\endlastfoot
\hline
field&purpose&field &purpose\\\hline
fuha&{\bf fu}nction {\bf ha}ndles, e.g., fuha.G, \ldots& nc&
{\bf n}umerical {\bf c}ontrols, e.g., nc.tol, \ldots\\
sw&{\bf sw}itches such as sw.bifcheck,\ldots&
sol&values/fields calculated at runtime\\
eqn&tensors $c,a,b$ for fast FEM setup&mesh&the geometry data and mesh\\
plot&switches and controls for plotting&file&switches etc for file output\\
time&timing information&pm&pmcont switches\\
fsol&switches for the {\tt fsolve} interface&
nu,np&\# PDE unknowns, \# meshpoints\\
u,tau&solution and tangent&branch&branch data\\\hline
usrlam&\multicolumn{3}{p{0.85\textwidth}|}{vector of user set target 
values for the primary parameter, default usrlam=[];} \\
mat&\multicolumn{3}{p{0.75\textwidth}|}
{problem matrices, in general data that is not saved to file, 
see Remark \ref{frefrem}}\\
\hline
\end{longtable}
\ece

\paragraph{Further general comments.} 
Concerning
the improved plotting, \pdepb~ uses telling axis labelling and, for
instance, a simplified user-friendly branch-plotting command: {\tt
  plotbra(p)}. By default, this plots the branch with
primary parameter on the $x$-axis and $L^2$-norm (now stored in
the internal part of the branch data) on the $y$-axis; the figure used
can be controlled by {\tt p.plot.brafig}. Similarly, {\tt
  plotbraf('p')} is now allowed for convenience and calls {\tt
  plotbra(p)} with structure {\tt p} from the file in directory {\tt
  'p'} with the highest point label. Moreover, in demo {\tt schnakfold} 
we provide some examples how to create movies of some continuation. 

Finally, we also added a wrapper 
to call \mlab's {\tt fsolve} routine; allthough this is typically slower 
than our own Newton loops, it may be useful, for instance, 
 to find solutions from poor 
initial guesses, see \S\ref{fchsec}. 

\paragraph{Acknowledgement} We thank Ben Schweizer (TU Dortmund) for help on the transformation to periodic boundary conditions used in \S\ref{pbcsec}.

\section{New features - by examples}
\subsection{Allen-Cahn model ({\tt acfold})}\label{s:ac}
As a first example we (re)consider the cubic--quintic 
Allen-Cahn equation from \cite[\S3.2]{p2p}, written as 
\huga{\label{ac}
-c\Delta u-\lam u-u^3+\ga u^5=0, 
}
on the rectangle $\Omega=(-1,1)\times(-0.9,0.9)$ with homogeneous 
Dirichlet \bcs. Our first task is to explain the new meaning of {\tt p.u}, 
parameter--switching and fold--continuation, and 
a new setup with a more efficient use of the FEM. The 
demo directory for this is {\tt acfold}. 

\subsubsection{Parameter switching}
 There are three parameters $c,\lam,\ga$, and 
in addition to the standard domain and \bcs \ setup known from \cite{p2p}, 
the init-routine {\tt acfold$\_$init} now initializes these and sets 
the primary continuation parameter:
{\small
\begin{verbatim}
% initialize auxiliary variables, here parameters of PDE
par(1)=1;         % linear cofficient of f
par(2)=0.25;      % diffusion coefficient
par(3)=1;         % quintic coefficient of f
p.u=[p.u; par'];  % augment p.u by parameters 
p.nc.ilam=1;      % set active parameter indices (here only one)
p.usrlam=[3.5 4]; % "target" values of the parameter
\end{verbatim}}
\noi
The functions defining $G$ and its Jacobian read 
{\small
\begin{verbatim}
function [c,a,f,b]=acfold_G(p,u)  % coefficient functions for AC
% separate pde and auxiliary variables, here "par", and interpolate to triangles
par=u(p.nu+1:end); u=pdeintrp(p.points,p.tria,u(1:p.nu)); 
c=par(2); a=0; b=0; f=par(1)*u+u.^3-par(3)*u.^5; end;
\end{verbatim}}
{\small
\begin{verbatim}
function [cj,aj,bj]=acfold_Gjac(p,u)  % jacobian for AC 
par=u(p.nu+1:end); u=pdeintrp(p.points,p.tria,u(1:p.nu));  
cj=par(2); bj=0; fu=par(1)+3*u.^2-par(3)*5*u.^4;  aj=-fu; end
\end{verbatim}}

\brem{\rm We recall, see \cite[Remark 3.2]{p2p}, that 
{\tt cj,aj,bj} in {\tt Gjac} are {\em not} the derivatives 
of $c,a,b$ in $G$. The notation only indicates that {\tt cj,aj,bj} 
are the coefficients needed to assemble $G_u$. In general, the 
relation between {\tt cj,aj,bj} and $c,a,b,f$ can be quite 
complicated, and only if $c,a,b$ are independent 
of $u$, and $f$ only depends on $u$ without derivatives 
(roughly: the semilinear case), 
then {\tt cj}$=c$, {\tt bj}$=b$, and {\tt aj}$=a-f_u$. Similar remarks 
apply to the functions {\tt p.fuha.bc} and {\tt p.fuha.bcjac}, 
see \S\ref{nlbcsec}. 
}\eex\erem

\subsubsection{Efficient use of FEM matrices in the semilinear case}\label{s:semi}
Exploiting a semilinear structure in the FEM assembling can give a
significant computational speedup: the FEM representation $G(u)=Ku-F$
of, e.g., $-\Delta u-f(u)$, can be obtained directly from 
$Ku$={\tt p.mat.K*u} and $F=${\tt p.mat.M*f(u)}, where {\tt p.mat.M} and {\tt
  p.mat.K} are the \textit{pre-assembled} mass and stiffness matrices, and {\tt f(u)} denotes
$f(u)$ as nodal values. In contrast, the FEM assembling via the
general routine {\tt [c,a,f,b]=G(p,u)}, calculates the coefficients
{\tt c,a,f,b} on the 
triangles after interpolation,
and then {\tt K}, {\tt F} are assembled from these at every Newton step.

In \pdepb~ the faster FEM setting is turned on by {\tt
  p.eqn.sfem=1}, which requires implementing the nodal routines for the 
Jacobian and residual, as well as setting the divergence tensor and,
if needed, the advection tensor. The matrices $M$ and $K$ are then 
generated via {\tt p=setfemops(p)} and stored in the structure {\tt p.mat}. 
For the {\tt acfold} demo the setup in {\tt acfold$\_$init} reads \\[2mm]
{\small 
{\tt p.sw.sfem=1; p.fuha.sG=@acfold$\_$sG; p.fuha.sGjac=@acfold$\_$sGjac;\\ 
p.eqn.c=1; p.eqn.b=0; p.eqn.a=0;} \\[2mm]
and the relevant routines are:
\begin{verbatim}
function r=acfold_sG(p,u)  
par=u(p.nu+1:end); u=u(1:p.nu); f=par(1)*u+u.^3-par(3)*u.^5;
r=par(2)*p.mat.K*u(1:p.nu)-p.mat.M*f; end 
\end{verbatim}\vs{-4mm}
\begin{verbatim}
function Gu=acfold_sGjac(p,u)  
par=u(p.nu+1:end); fu=par(1)+3*u.^2-par(3)*5*u.^4; Fu=spdiags(fu,0,p.nu,p.nu);
Gu=par(2)*p.mat.K-p.mat.M*Fu; end
\end{verbatim}
}
For problems involving the advection tensor $b$, analogously define 
{\tt p.eqn.b} and use the matrix {\tt p.mat.Kadv} in the routines; 
see the {\tt acfront} and {\tt acffold} demos and, for a system, the {\tt schnaktravel} and {\tt nlb} demos. 
Assembling the mass matrix at startup, and automatic updates at mesh adaption or refinement, are controlled by setting {\tt p.sw.sfem} to a nonzero value. If the user only wants to use pre-assembled mass matrix (hence no fast FEM for \eqref{gform}), this would be {\tt p.sw.sfem=-1}. 

\brem\label{r:cust}(Customisation) {\rm If the operators in \eqref{gform} depend on parameters, the ``semilinear'' implementation explained above can sometimes be extended by splitting the operators suitably. For instance, if the operator $c$ can be written as $c=c_1 + \lambda c_2$, the routine {\tt setfemops} can be locally (in the demo directory) modified to generate {\tt p.mat.K} and {\tt p.mat.K2} so that the residual reads {\tt p.mat.K+lam*p.mat.K2}. See \S\ref{s:vkplate} for an example.}
\eex\erem

\brem\label{frefrem}(Mesh refinement, saving of FEM matrices) {\rm 
If a (semilinear) 
problem is to be run on a fixed mesh, then setting {\tt p.sw.sfem=1} and 
setting {\tt p.fuha.sG} and {\tt p.fuha.sjac} as above can replace 
the old setting with {\tt p.fuha.G} completely. However, if adaptive 
mesh refinement is desired, then {\tt p.fuha.G} is still needed to 
identify the triangles to be refined. The required new matrices 
{\tt p.mat.M}, {\tt p.mat.K} and {\tt p.mat.Kadv} are automatically 
reassembled during mesh adaption. Moreover, 
to save hard disk space, in the standard setting 
{\tt p.fuha.savefu=@standsavefu}, the struct {\tt p.mat} is not saved 
in the solution files. When loading via {\tt loadp}, it is automatically 
regenerated. However, when loading a {\tt p} struct into the \mlab~ 
workspace by a double click, this is not the case and a manual call 
to {\tt setfemops} is needed.
}
\eex\erem 

\brem\label{tintrem1}(time integration) {\rm Although \pdepb\ is not primarily intended for 
time integration, we also 
provide some extensions of the simple general 
time--integrator {\tt tint} from \cite{p2p} for systems $\pa_t u=-G(u,\lam)$. 
{\tt tintx} uses the same (linearly implicit) algorithm 
as {\tt tint}, based on the full {\tt p.fuha.G} syntax, 
but also returns a time--series of the residual at each plotting 
step, and saves the time evolution in {\tt p.file.pre}. In detail, 
at startup, the full structure 
{\tt p} is saved (if it does not yet exist) to ``pre/pt0.mat'', 
and afterwards only {\tt p.u} is saved at the selected time steps. 
To load a point {\tt ptn} from ``pre'', we then use a modified 
{\tt p=loadp2(pre,'ptn','pt0')}. 

For the fast FEM setting {\tt sfem=1} the integrators {\tt tints} and {\tt tintxs} are much more efficient -- typically by at least a factor 10. 
Again these are simple linearly implicit schemes, but based on an LU 
decompostion of $M+dt K$. See \S\ref{tintsec} for more comments. }\eex
\erem 

\subsubsection{Fold detection, point types and parameter switching}\label{s:parsw}
\begin{figure}[!h]
\bce
\begin{tabular}{p{0.22\textwidth}p{0.22\textwidth}p{0.22\textwidth}p{0.22\textwidth}}
(a)&(b)&(c)&(d)\\
\ig[height=36mm]{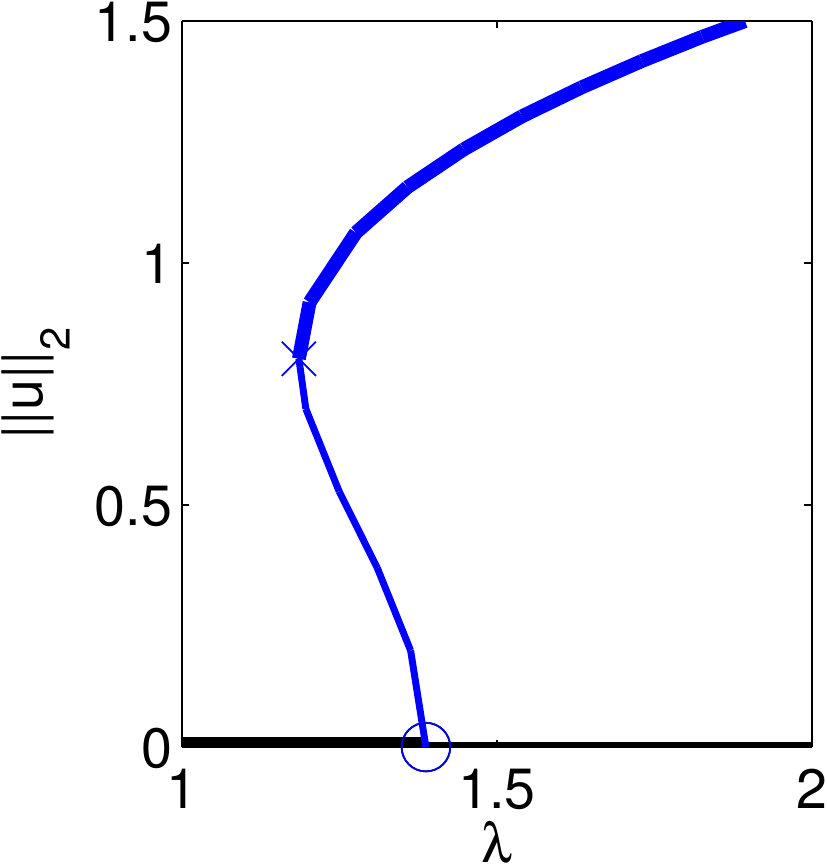}&\ig[height=36mm]{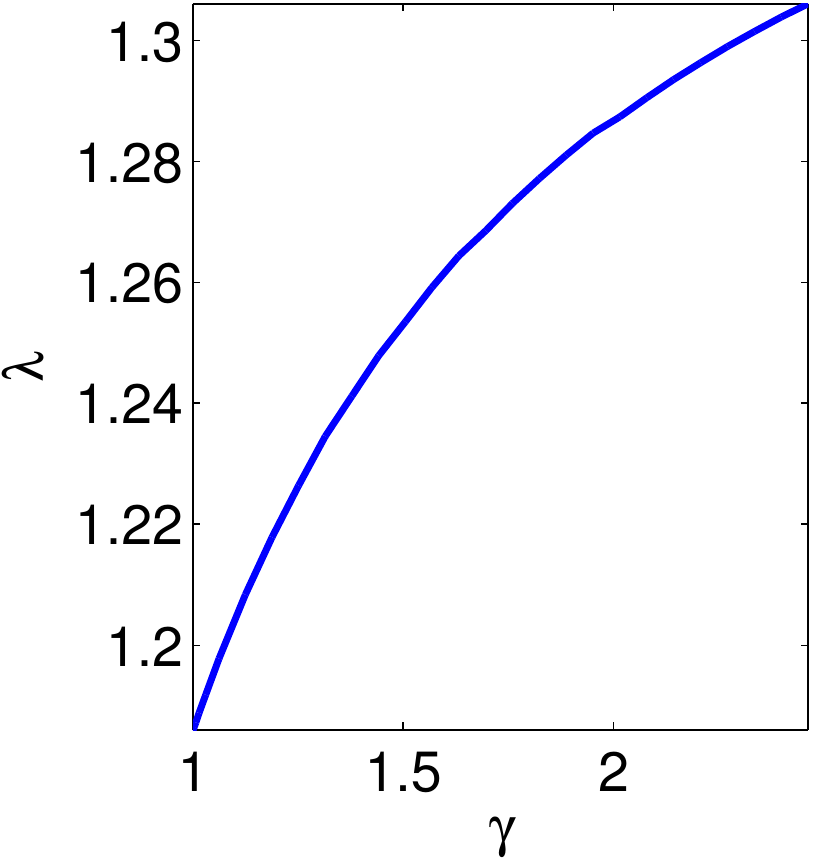}&
\raisebox{2mm}{\ig[height=35mm]{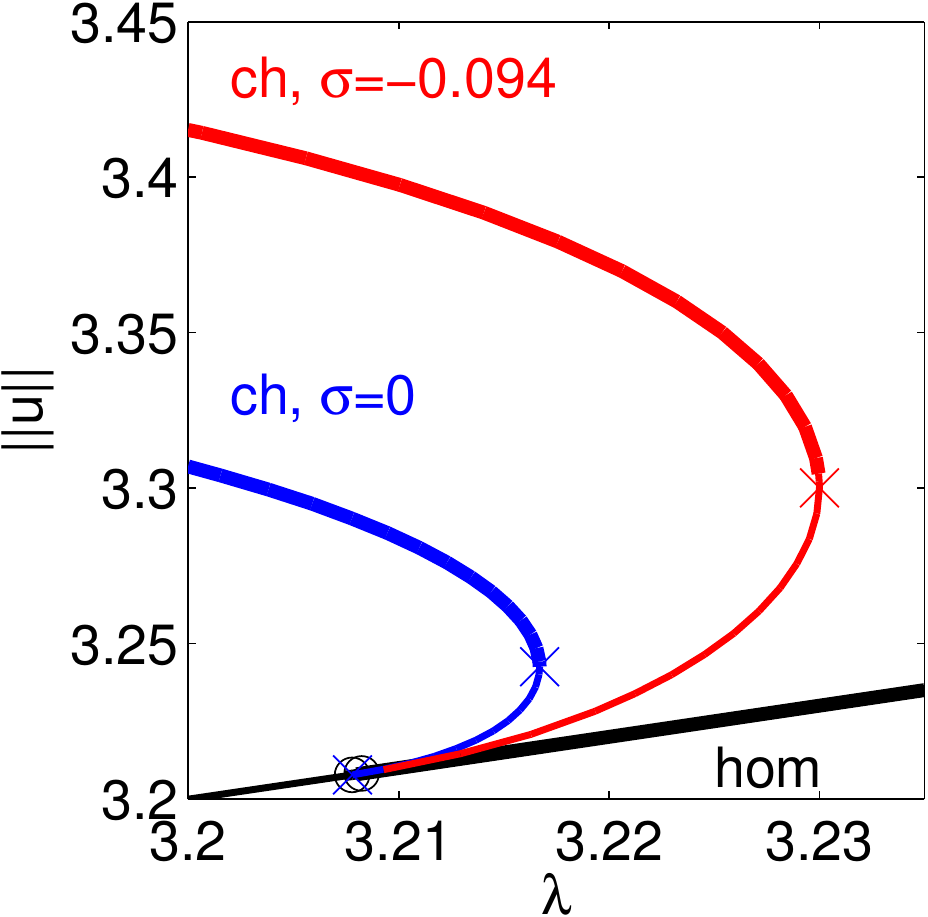}}&
\raisebox{2mm}{\ig[height=35mm]{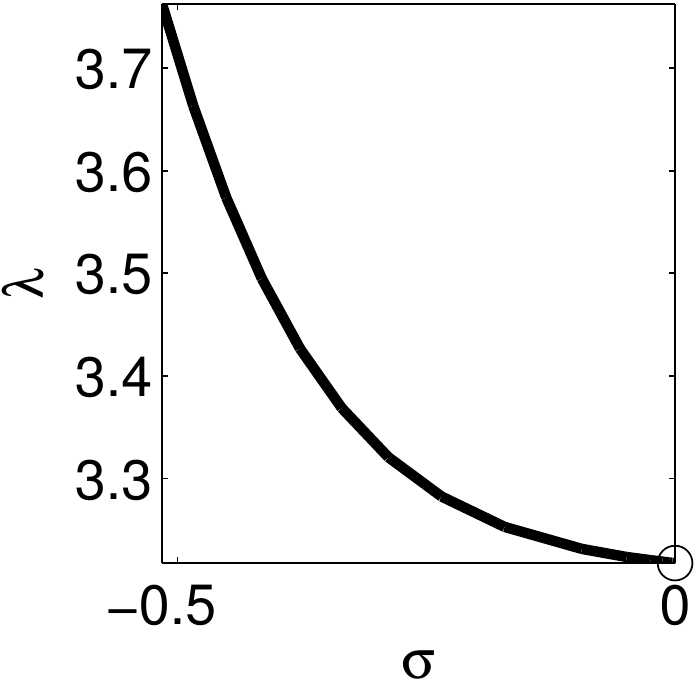}}
\end{tabular}
\ece

\vs{-5mm}
\caption{{\small Fold continuation in the  Allen--Cahn model \reff{ac} and 
the Schnakenberg model \reff{mod1}. 
(a),(b) First bifurcating branch ($(c,\ga)=(0.25,1)$) 
and ``fold position'' $\lam$ over the 
quintic parameter $\ga$ in \reff{ac}. 
(c) ``cold hexagon'' solution branch (blue) for \reff{mod1} with 
$\sig{=}0$. 
(d) continuation of the fold for $\sig{=}0$ from (c) in $\sig$; 
afterwards, the red branch in (c) 
was obtained via {\tt foldexit} at $\sig\approx-0.094$ 
and continuing in $\lam$ again, with positive and negative {\tt ds}. 
For details see {\tt acfold\_cmds.m} and {\tt schnakcmds.m}, respectively.}
\label{sffig}}
\end{figure}

Coming back to \eqref{ac}, after locating the well-known bifurcation points 
(eigenvalues of the Dirichlet Laplacian)
from the trivial branch $u\equiv 0$, we switch in {\tt acfold$\_$cmds} to the first bifurcating branch and continue it including fold-detection by 
{\small 
\begin{verbatim}
q=swibra('p','bpt1','q',-0.2); q.sw.foldcheck=1; p=cont(p);
\end{verbatim}
}
\noindent
where fold detection works by bisection as for branch points. 
The resulting branch is plotted in Figure~\ref{sffig}(a) with the fold 
point marked. It is also stored in the file {\tt q/fpt1.mat} and assigned 
a special point type in the branch {\tt p.branch}. 
These point types in \pdepb~ branches are:

\begin{tabular}{rl}
      -1&= initial point or restart\\
      -2 &= guess from {\tt swibra} for the initialization of branch switching\\
       0 &= regular point\\
       1 &= bifurcation point (found with {\tt bifdetec})\\
       2 &= fold point (found with {\tt folddetec})
\end{tabular}

\medskip
Switching parameters in order to continue a stored solution in the previously constant diffusion rate $c$ (parameter number 2 in the implementation) goes simply by {\tt w=swiparf('q','pt10','w',2);} where the essential change done by {\tt swiparf} is setting {\tt w.ilam=2} and where {\tt w} is the name of the new branch and {\tt q/pt10} is the file name of the stored solution. Before continuation by {\tt w=cont(w);} some adjustments to the settings are useful in this case:
{\small
\begin{verbatim}
w.nc.lammin=0.1; w.sol.ds=-0.01; w.sol.xi=1e-6; 
\end{verbatim}
}
\noindent where the small weight $\xi$ is useful since the problem is more sensitive in the diffusion coefficient.

\subsubsection{Fold and branch point continuation}\label{s:acfold}

We explain continuation of the fold-point in the {\tt q}--branch in 
Fig.~\ref{sffig}(a). Constraining continuation to folds requires an 
additional free parameter, e.g, $w:=(\lam,\ga)$. Before going into 
the practice in \pdepb\ we briefly discuss the background of 
fold--(and branchpoint--) continuation. In this case \pdepb\ discretizes 
the extended system  
\huga{\label{heq} 
H(U)=\bpm G(u,w)\\ \partial_u G(u,w)\phi\\ 
\|\phi\|^2_{L^2}-1\\
p(U)\epm=0, \quad U=(u,\phi,w),
}
so that $\phi$ is in the kernel of $\partial_u G$ with 
$L^2$-norm constrained to 1 by the third equation, 
and $p(U)=0$ is the arclength equation \eqref{e:arc}. Thus 
the FEM discretization of \reff{heq} is a system of {\tt p.nu+p.nu+2} equations 
in {\tt p.nu+p.nu+1} unknowns. 

For continuation of \reff{heq} we need the Jacobian 
\huga{\label{spjac} 
D_U H(U)=\bpm \pa_u G& 0&\pa_w G\\
\pa_u(\pa_u G\phi)&\pa_u G&\pa_w (\pa_u G\phi)\\
0 & 2\phi^T & 0\\
\xi \dot{u}^T&\xi \dot{\phi}^T&(1-(\xi+\xi_q)/2)\dot{w}\epm,
}
where depending on {\tt p.sw.jac} and {\tt p.sw.sfem} 
$\pa_u G$ is calculated numerically or 
assembled using {\tt p.fuha.Gjac} or {\tt p.fuha.Gjac}, 
respectively, $\phi$ only occurs linearly 
in \reff{heq}, derivatives with respect to $w$ are done via finite differences, and the computationally most costly part is the 
evaluation of $\pa_u(\pa_u G\phi)$. While this is done numerically for 
{\tt p.sw.spjac=0}, the user is urged to implement $\pa_u(\pa_u G\phi)$ in 
a routine {\tt p.fuha.spjac} and set {\tt p.sw.spjac=1}. 

\medskip
In the {\tt acfold} demo, this is readily done since $\pa_u(\pa_u G\phi)=f_{uu}\phi$ so that we can use a pre-assembled mass matrix {\tt p.mat.M} as in the 
FEM assembling for semilinear problems discussed in \S\ref{s:semi}:
{\small
\begin{verbatim}
function Guuph=acfold_spjac(p,u) 
ph=u(p.nu+1:2*p.nu); par=u(2*p.nu+1:length(u)); u=u(1:p.nu); 
fuu=6*u-20*par(3)*u.^3; Guuph=-(p.mat.M*diag(fuu))*diag(ph); 
\end{verbatim}
}
The use of the extended system \eqref{heq} and its Jacobian for subsequent continuations (the fold-/branch point-continuation mode) in \pdepb~ is turned on by calling {\tt spcontini}; for instance (see {\tt acfold$\_$cmds}):
{\small\begin{verbatim} 
qf=spcontini('q','fpt1',3,'qf');% init fold continuation with par 3 as new active parameter
qf.plot.bpcmp=3; clf(2);        % use this new parameter for plotting
qf.nc.tol=1e-5;                 % increase tolerance as typically required for fold cont.
qf.sol.ds=1e-3;                 % new stepsize in new primary parameter
\end{verbatim}}
The branch computed now by  {\tt qf=cont(qf)} is plotted in 
Fig.~\ref{sffig}(a). Note here {\tt p.nc.ilam=[1;3]}. Starting a normal continuation from a point stored during a fold- or 
branch point-continuation is done on by calling {\tt spcontexit} as in:
{\small\begin{verbatim} 
q1=spcontexit('qf','pt10','q1'); q1.nc.tol=1e-8; q1.sol.ds=1e-3; q1=cont(q1);
\end{verbatim}}

\brem\label{rem:spcont}{\rm 
The switch {\tt sw.spcont} is internally used to distinguish these modes: 0 means normal, 1 means branch point and 2 means fold point continuation, in agreement with the point types mentioned above. This stores the continuation mode since during continuation of fold/branch points the point types are normally zero as for usual continuation. For branch switching from a branch point continuation it can be useful to generate a new guess for the tangent vector {\tt tau}. This can be conveniently done with the new routine {\tt p=getinitau(p)} and the new method to call {\tt swibra} as in {\tt q=swibra(p,0.05,'q')}. See the demo {\tt bratu} for an example.}
\eex\erem

\subsection{Semilinear structure in a system: the Schnakenberg model ({\tt schnackfold})}\label{ssec}
\subsubsection{Fold continuation}
As an example of a system (i.e. $N=${\tt neq}$>1$) we consider the Schakenberg model
\huga{\label{mod1}\pa_t U=
D\Delta U+N(u,\lam)+\sigma \left(u-\frac 1 v\right)^2\bpm 1\\-1\epm, 
\quad N(U,\lam)=\bpm-u+u^2v \\\lambda -u^2v\epm, 
}
with $U=(u,v)(t,x,y)\in\R^2$, diffusion matrix 
$D=\ssmatrix{1 & 0\\ 0 & d}$, 
$d$ fixed to $d=60$, and  bifurcation parameters $\lambda\in \R_{+}$ 
and $\sig\in\R$. 
System \reff{mod1} has the homogeneous stationary solution $(u,v)=(\lam,1/\lam)$, 
which becomes Turing unstable for $\lam\le \lam_c\approx 3.2085$, 
independent of $\sig$, with critical wave-vectors $k=(k_1,k_2,k_3)$ with 
$|k|=k_c=\sqrt{\sqrt{2}{-}1}{\approx}0.6436$. 
Here $\sig$ can be used to turn certain 2D bifurcations from sub--to supercritical, and many branches of patterns exhibit one or many folds (``snaking'') \cite[\S4.2]{p2p} and \cite{uwsnak14}. 
Fold continuation can be used to discuss snaking widths, see 
\cite{uwfc14}, which however, requires rather large systems with grid size $\CO(10^5)$ so that finite differences for $\pa_u(\pa_u G \phi)$ 
in \eqref{spjac} are inefficient. 

Following the approach discussed in \S\ref{s:semi} for semilinear problems, here we implement $\pa_u(\pa_u G \phi)$ in the simplified nodal FEM format (in the demo {\tt schnakfold} also the PDE itself is implemented according to \S\ref{s:semi}). Denoting 
 $u=(u_1,u_2)$, $\phi=(\phi_1,\phi_2)$, $f=(f_1,f_2)$ for the components of $u,\phi$ and $f$, we have 
\huga{
\pa_u(\pa_u G \phi)=\bpm (\pa_{u_1}^2f_1)\phi_1+(\pa_{u_1}\pa_{u_2}f_1)\phi_2& 
(\pa_{u_1}\pa_{u_2}f_1)\phi_1+(\pa_{u_2}^2f_1)\phi_2\\
(\pa_{u_1}^2f_2)\phi_1+(\pa_{u_1}\pa_{u_2}f_2)\phi_2& 
(\pa_{u_1}\pa_{u_2}f_2)\phi_1+(\pa_{u_2}^2f_2)\phi_2
\epm. 
}
Using the nodal values for $\pa_i\pa_jf_k$ and multiplication with the mass matrix {\tt p.mat.M}, 
this is implemented in {\tt schnakspjac.m} of the {\tt schnakfold} demo. The continuation in \pdepb~ works as discussed in \S\ref{s:acfold}. 
We plot the ``cold hexagon'' branch in $\sig$ on a small domain and the 
continuation of its first fold point in Figure~\ref{sffig}(c),(d).  

\subsubsection{Time integration and movies}\label{tintsec}
As indicated in Remark \ref{tintrem1}, we also use 
the Schnakenberg model as an example for time integration with 
{\tt tints} and {\tt tintxs}, see {\tt schnakcmds2.m}. 
In {\tt p=tints(p,dt,nt,pmod,} {\tt nffu,varargin)}, for time--integration with 
stepsize ${\tt dt}=\del$ of 
the FEM representation $Mu_t=Ku-Mf(u)$ with a $u$-independent $K$ we use 
$$
\Lam u^{n+1}=\left[Mu^{n}+\del Mf(u)\right], \text{ where }
\Lam=M+\del K, 
$$
To solve this linear system for $u^{n+1}$, we $LU$-decompose $\Lam$ at 
startup. The nodal values $f(u)$, which are also needed in 
{\tt p.fuha.sG}, must be encoded in {\tt f=nffu(p,u)}, 
$M$ is taken from {\tt p.mat.M}, and for 
$K$ there are the following options: If {\tt varargin=[]}, then 
$K$ (including advective terms, if non--zero) 
is built from {\tt p.mat.K} and {\tt p.mat.Kadv}, i.e., 
$K={\tt p.mat.K}+{\tt p.mat.Kadv}$. 
If {\tt varargin=K}, then $K={\tt K}+{\tt p.mat.Kadv}$, 
and if {\tt varargin=[K,Kadv]}, then $K={\tt K}+{\tt Kadv}$. 
This is useful if diffusion parameters are not included in $K$ 
but used explicitly in {\tt p.fuha.sG}, see {\tt acfold} for 
an example. 

\begin{figure}[!h]
\bce
\begin{tabular}{p{30mm}p{70mm}}
(a)&(b)\\
\ig[height=41mm]{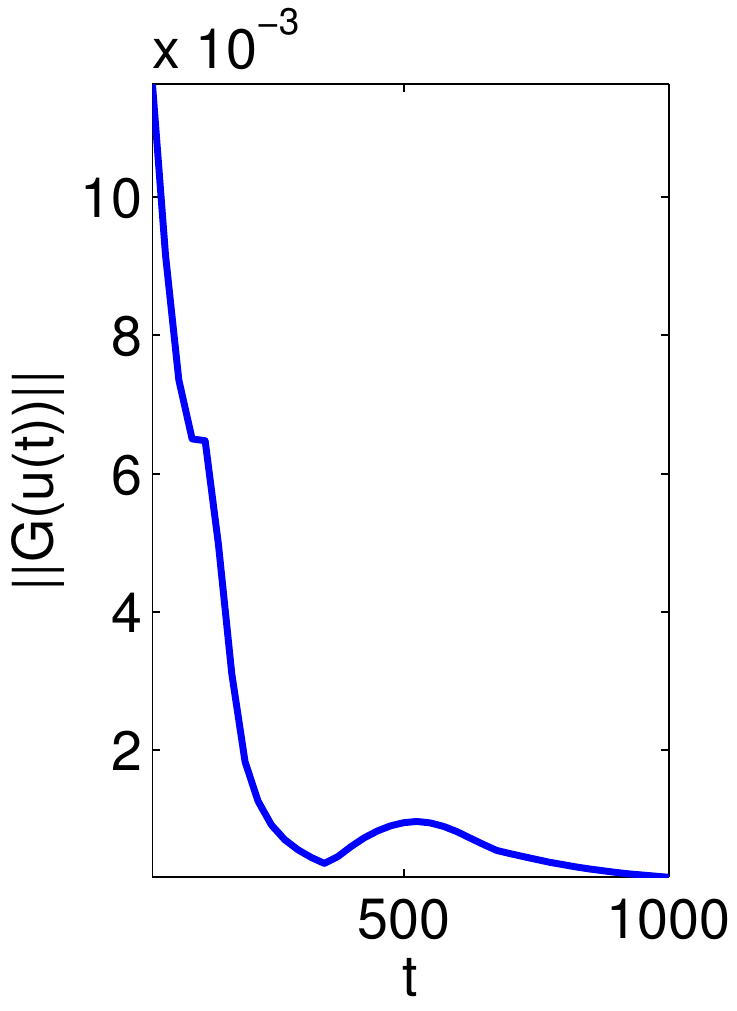}&
\raisebox{22mm}{\begin{minipage}{80mm}\ig[width=70mm,height=23mm]{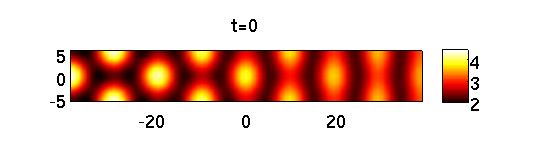}
\\[-4mm]
\ig[width=70mm,height=25mm]{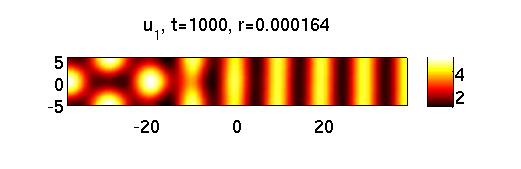}\end{minipage}}
\end{tabular}
\ece

\vs{-5mm}
\caption{{\small Time integration of \reff{mod1}, $\lam=2.8,\sig=0$. 
(a) time evolution of the residual; (b) 
initial guess for a front between hexagons and stripes in \reff{mod1}, 
and the solution at $t=1000$. The latter can be used as starting point for 
continuation of the stationary problem in, e.g., $\lam$, which gives 
a ``snaking'' branch of stationary fronts. See \cite{uwsnak14} for 
further discussion, including the bifurcation of such snaking branches 
of localized solutions from so--called bean branches. In (b) the horizontal and vertical axes are $x=x_1$ and $y=x_2$ resp. \label{stfig}}}
\end{figure}

If applicable, {\tt tints} is at least 10 times faster than the 
old {\tt tint} method which assembles $M,K,F$ at each time-step, 
and cannot take advantage of a precomputed $LU$-decomposition of $\Lambda$.
Besides the simple--interface versions {\tt p=tint(p,dt,nt,pmod)} and 
{\tt p=tints(p,dt,nt,pmod,nnfu,varargin)} 
there are also versions {\tt  tintx} and 
{\tt  tintxs} which write the solution at selected time--steps in a file, 
and return some diagnostics, such as the time series of the 
residual $\|G(u(t))\|$. See Fig.~\ref{stfig} for an example. 
Again we remark that all these time--integrators should be seen as 
templates for more problem-adapted routines; in particular, there 
is no error or stepsize control. 

Finally, in {\tt moviescript.m} 
we give some examples for movie creation. Typically, these require some 
customized {\tt plotsol} and {\tt plotbra}, see {\tt mplotsol.m} 
and {\tt mplotbra.m}, but otherwise explain themselves. 
See also the end of {\tt schnakcmds2.m} for time-integration movies. 

\subsection{Nonlinear boundary conditions ({\tt nlbc})}\label{nlbcsec}
We consider 
\huga{\label{nlbc}
-\Delta u=0\text{ in }\Om,\quad \pa_n u+\lam s(x,y) f(u)=0,
}
taken from \cite{MA11}, where  
we choose $s(x,y)=0.5+x+y$, $f(u)=u(1-u)$, 
and, for diversity, $\Om=D=\{(x,y)\in\R^2:\|(x,y)\|_2<1\}$. Thus we have the simple 
linear Laplace equation with nonlinear boundary conditions, where 
we take $\lam$ as our bifurcation parameter. 
Clearly, $u\equiv 0$ and 
$u\equiv 1$ are two trivial branches, $f(u)>0$ in between, and a 
crucial feature is that 
the weight function $s$ changes sign on $\pa \Om$. Moreover, 
$\spr{s}:=\int_{\pa\Om} sd\Gamma>0$, which corresponds to the case $\spr{s}<0$ 
in \cite{MA11}; we chose the different sign to make the connection to 
the general form \reff{gnbc}, 
${\bf n}\cdot (c \otimes\nabla u) + q u = g$, of boundary conditions 
more transparent, i.e., we choose  $q=\lam s(x,y) (1-u)$ and $g=0,c=1$. 

As the model is related to gene frequencies, in applications 
one is mostly interested in solutions $u$ with $0\le u(x,y)\le 1$, and in 
\cite{MA11} a number of remarkable results are shown, essentially  
assuming $f$ of the above form and $\spr{s}\ne 0$. In particular, 
if $\spr{s}>0$ (in our convention), the only 
nontrivial solutions with $0\le u(x,y)\le 1$ in $\Om$ are on a global branch 
$u(\cdot;\lam)$, $\lam> \lam_0>0$, and these are exponentially stable in the heat equation associated with \eqref{nlbc}.  
Additionally, for $0<\lam\le\lam_0$, the trivial solution $u\equiv 0$ 
is stable. 
We do recover these results numerically, but for completeness we 
drop the restrictions $\lam>0$ and  $0\le u(x,y)\le 1$; note, for 
instance, that at $\lam=0$ any constant $u$ is a solution of \reff{nlbc}. 

\begin{figure}[!h]
{\small 
\begin{tabular}[t]{p{40mm}p{120mm}}
\ig[width=39mm,height=40mm]{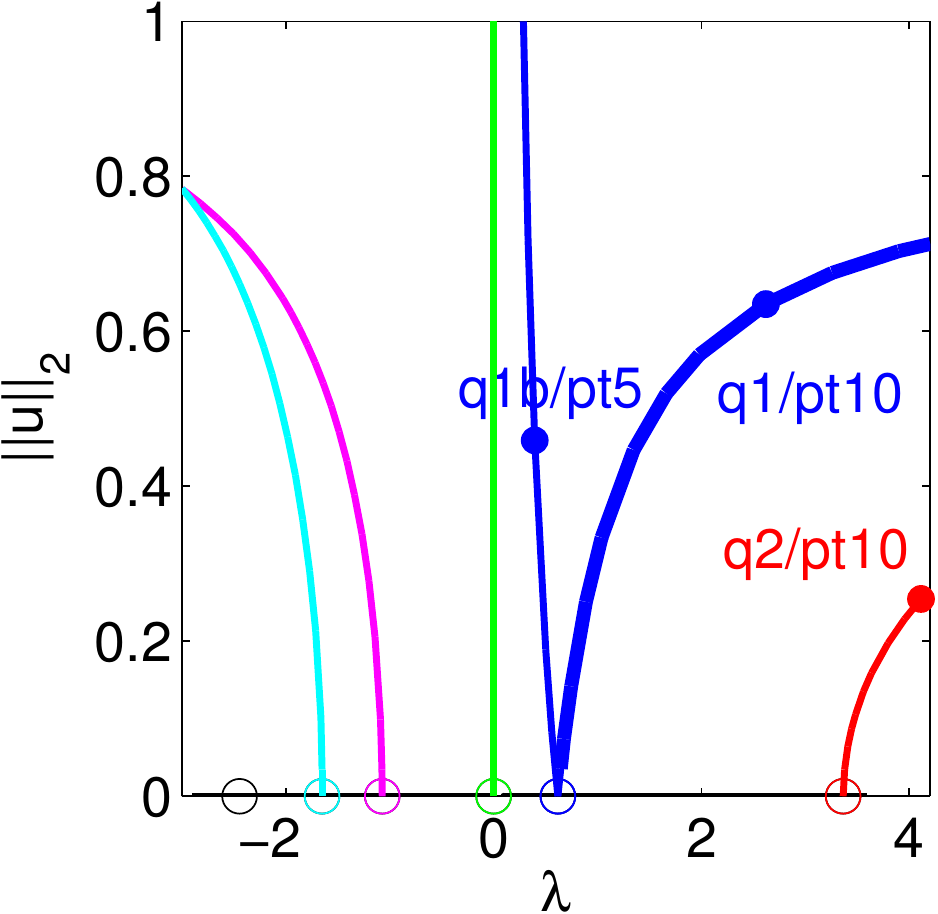}&\ig[width=38mm]{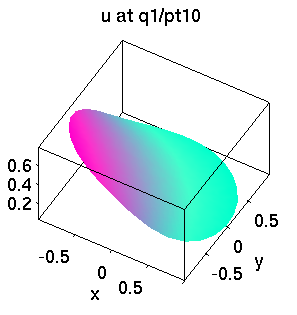}
\ig[width=38mm]{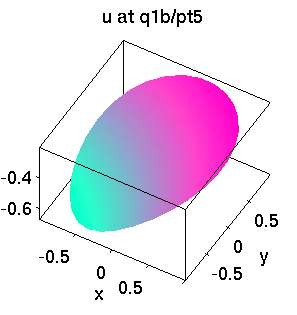}
\ig[width=38mm]{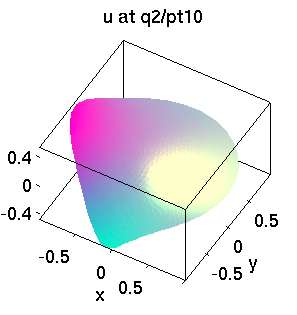}
\end{tabular}
}
\caption{{\small Bifurcation diagram (BD) and example plots for \reff{nlbc}. 
Here at $\lam_0\approx 0.62$ the bifurcation is transcritical, 
but $0\le u\le 1$ holds  only on the blue branch to the right of $\lam_0$, i.e. branch {\tt q1}. 
$u\equiv 0$ is stable for $0\le\lam\le\lam_0$. At $\lam=0$ a 
vertical branch of constant solutions bifurcates, and there 
are further bifurcation points both for $\lam<0$ and $\lam>\lam_0$, 
but the {\em only} non--constant solutions $u$ with $0\le u\le 1$ 
are on {\tt q1}. 
\label{nlbdf1}}}
\end{figure}

The coefficients {\tt c=1; a=0; b=0;} in {\tt fuha.G} and in {\tt fuha.Gjac}
are rather obvious, and it remains to encode the boundary conditions in the 
form \reff{gnbc}. The format of the ``boundary matrix'' {\tt bc} 
in the \ptool\ is rather 
unhandy, which is why we provide the functions {\tt bc=gnbc(neq,varargin)} 
and {\tt bc=gnbcs(neq,varargin)}, see \cite[\S3.1.4]{p2p}. Here we need 
the second version which takes string arguments containing expressions 
in {\tt x,u}, and set up {\tt fuha.bc} as 
{\small\begin{verbatim} 
function bc=nlbc(p,u) % nonlin., x-dep. BC; 
lam=u(p.nu+1); enum=max(p.mesh.e(5,:)); % find number of edges 
g=mat2str(0);q=[mat2str(lam) '*(0.5+x+y).*(1-u)']; bc=gnbcs(p.nc.neq,enum,q,g); 
\end{verbatim}}
The function {\tt fuha.bcjac} must provide the coefficients to assemble 
the $u$ derivatives of the \bcs. Accordingly, 
{\small\begin{verbatim} 
function bc=nlbcjac(p,u) % generate bc-matrix for derivatives of BC 
lam=u(p.nu+1); enum=max(p.mesh.e(5,:));
g=mat2str(0);qj=[mat2str(lam) '*(0.5+x+y).*(1-2*u)']; bc=gnbcs(p.nc.neq,enum,qj,g); 
\end{verbatim}}
With these definitions, \reff{nlbc} can now be run in \pdepb\ with 
{\tt sw.jac=1} (assembled Jacobians) in a standard way, see {\tt 
nlbccmds.m}, and Fig.~\ref{nlbdf1} for some results. 

\brem\label{nlbcrem}{\rm  
There is less flexibility when using 
{\tt fuha.bcjac} for linearizing \bcs\ than in using {\tt fuha.Gjac} for 
linearizing $G$, as $c$ in \reff{gnbc} is fixed by \reff{gform}. 
Thus, essentially {\tt bcjac} can be used if only 
$q$ or $g$ in \reff{gnbc} depend on $u$, and 
in more general cases one has to use {\tt sw.jac=0}, i.e. numerical differentiation, where 
{\tt fuha.bcjac} is not used. On the other hand, in the case of 
linear homogeneous \bcs, one has {\tt bcjac=bc} and hence 
can set {\tt fuha.bcjac=fuha.bc}. This is what we do in most examples. 
}
\eex\erem 

\subsection{Integral constraints: the functionalized Cahn-Hilliard equation ({\tt fCH})}\label{fchsec}
\begin{figure}[!h]
{\small 
\begin{tabular}[t]{p{42mm}p{110mm}}
(a) BD&(b) \\
\ig[width=39mm]{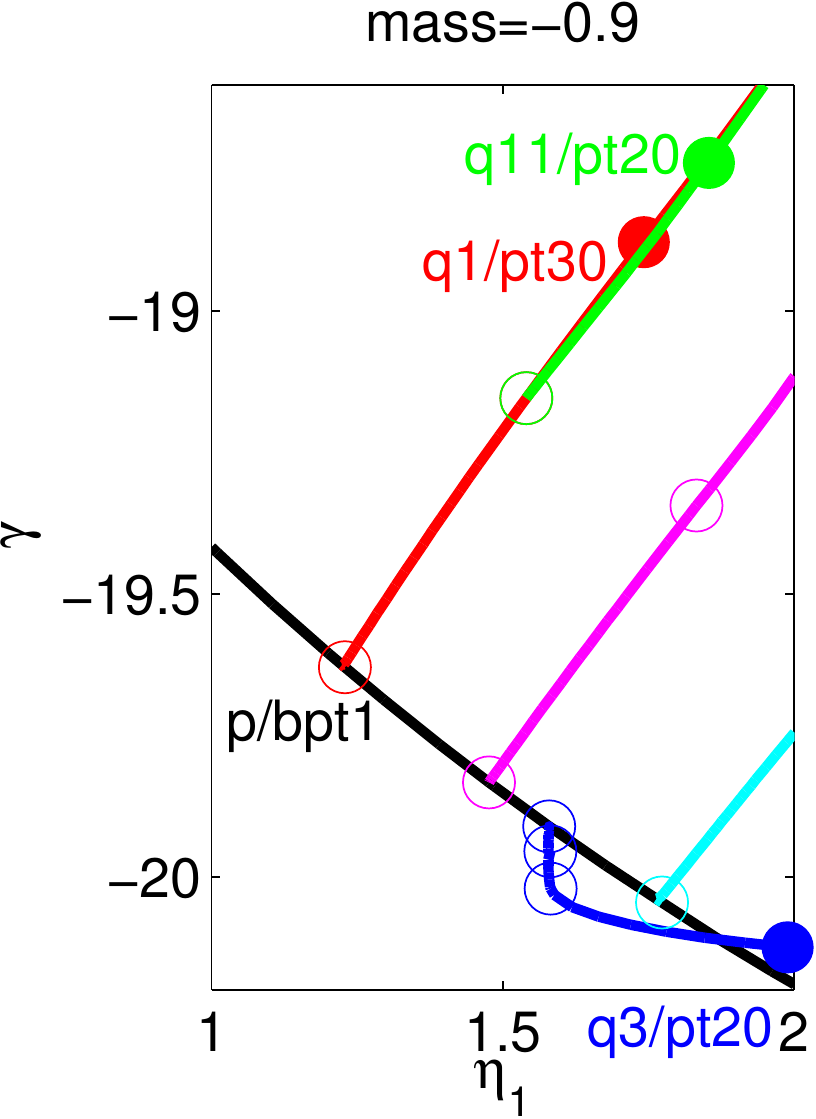}&\ig[width=33mm,height=46mm]{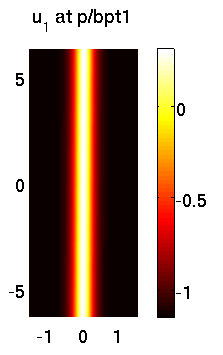}
\ig[width=21mm,height=46mm]{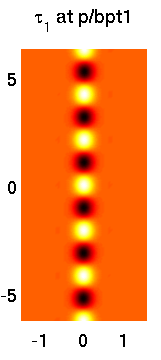}
\ig[width=21mm,height=46mm]{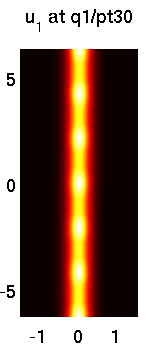}
\ig[width=21mm,height=46mm]{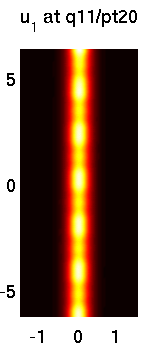}\\
(c)&(d)\\
\ig[width=20mm,height=46mm]{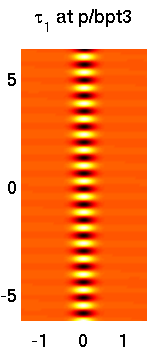}
\ig[width=20mm,height=46mm]{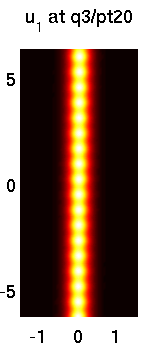}&
\ig[width=48mm]{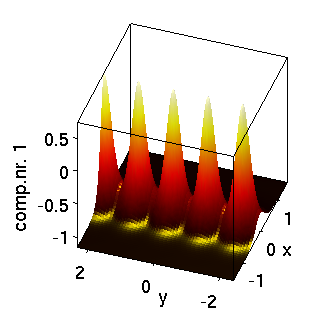}
\end{tabular}
}
\caption{{\small Bifurcations from a straight channel ({\tt fchcmds1.m}). 
(a) BD for 
 $\eps=0.35, \eta_2=2.5$, $m=-0.9$. (b) $u_1$ at 1st BP, tangent direction, 
pearled branch, and solution on secondary branch. (c) tangent a 3rd BP, and 
solution on bifurcating branch. Colormap everywhere 
as in the first plot in (b). In (b) and (c) the horizontal and vertical axes are $x=x_1$ and $y=x_2$ resp. For (d) see text.} \label{fchf1}}
\end{figure}

\begin{figure}
\begin{tabular}{ll}
\ig[height=40mm]{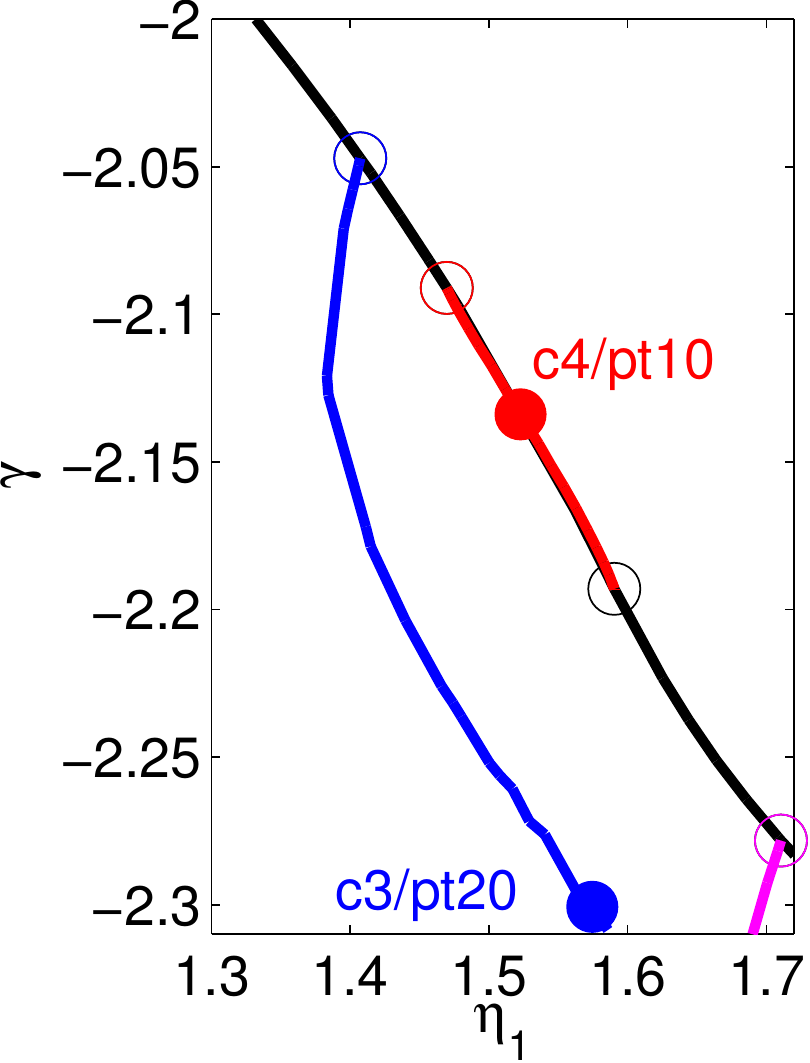}
\ig[height=40mm]{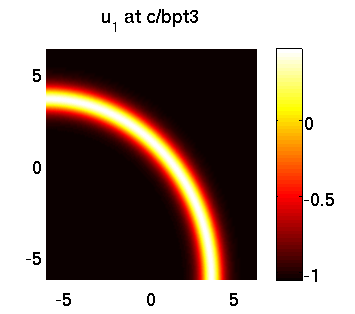}
&\ig[height=40mm]{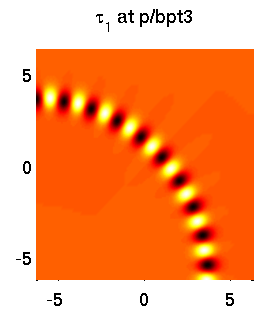}
\ig[height=40mm]{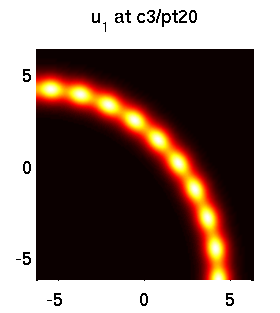}\\
\ig[height=40mm]{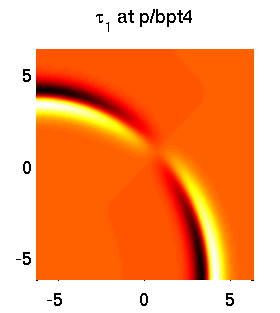}
\ig[height=40mm]{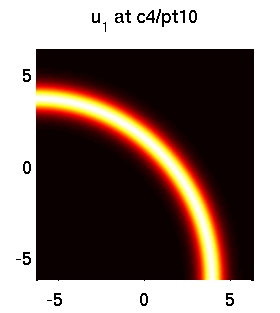}
&\ig[height=40mm]{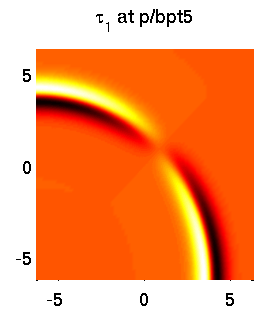}
\ig[height=40mm]{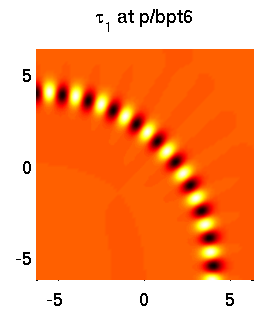}
\end{tabular}
\caption{{\small Bifurcation diagram and example directions/solutions 
for a curved channel. $\eps=0.75, \eta_2=2$, and $\eta_1=1$ initially, 
with a quarter-circle channel initial guess with mass $m\approx -0.8675$.  
On continuation in $\eta_1$, the blue branch bifurcates at the 3rd bifurcation point (BP); 
the first two BPs also yield pearling. The 4th and 5th BPs are connected 
by a ``meandering'' branch, and the 6th BP yields pearling again. The horizontal and vertical axes in the solution plots are $x=x_1$ and $y=x_2$ resp.
See {\tt fchcmds2.m}.} \label{fchf1c}}
\end{figure}

\begin{figure}
{\small 
\begin{tabular}[t]{p{110mm}p{35mm}}(a)&(b)\\
\ig[height=40mm,width=30mm]{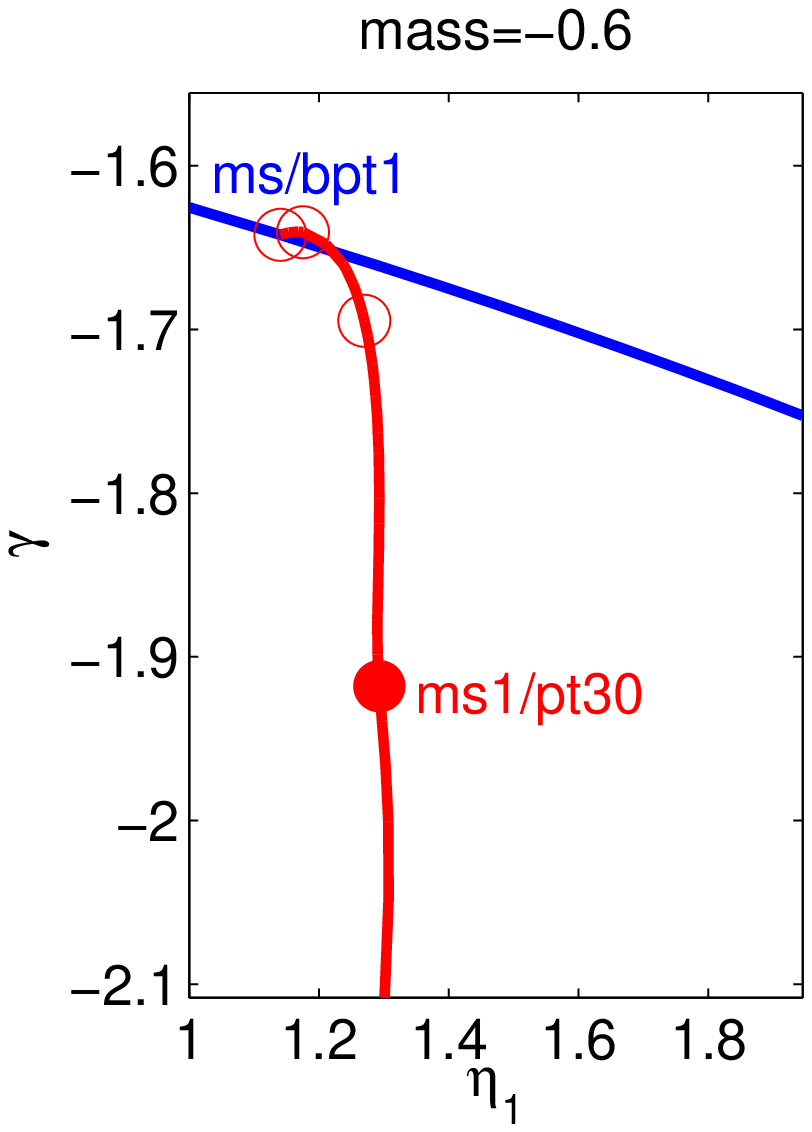}
\ig[height=40mm,width=28mm]{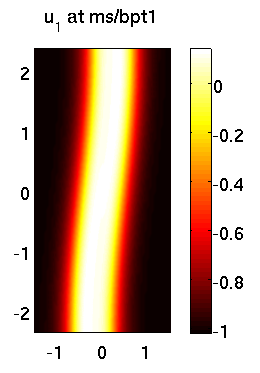}
\ig[height=40mm,width=20mm]{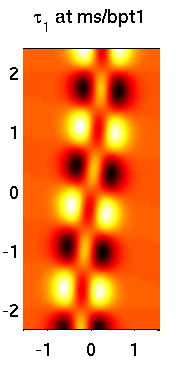}
\ig[height=40mm,width=20mm]{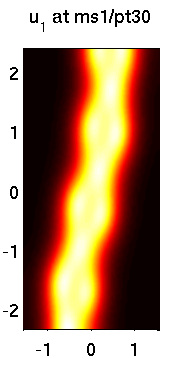}
&\ig[height=40mm]{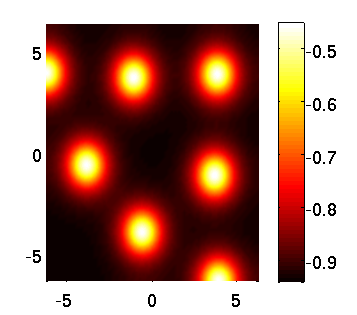}
\end{tabular}
}
\caption{(a) Bifurcation from a ``thick'' meandering branch. 
$(\eps,\eta_3,m)=(0.35,3,-0.6)$, 
initially $\eta_1=1$ and initial guess of type \reff{fchiguess} with 
$a=1.25, a_2=0.3, b_1=0.5, b_2=10$. Note that the tangent at bifurcation 
is different from the previous pearling in the sense that it also 
has a transversal periodic structure. 
(b) $u$ for spots obtained 
from a ``sharp bend'' initial guess, $(\eta_1,\eta_2,\eps)=(1,2,0.8)$, 
$m\approx -0.85$. See {\tt fchcmds3.m} for details. The horizontal and vertical axes in the solution plots are $x=x_1$ and $y=x_2$ resp.
  \label{fchf1b}}
\end{figure}

As an example of a problem with a constraint we consider 
the so called functionalized Cahn--Hilliard equation from \cite{DHPW12}, 
\huga{\label{tfch} 
\pa_t u=-\CG[(\eps^2\Del-W''(u)+\eps\eta_1)(\eps^2\Del u-W'(u)+\eps\eta_dW'(u))], 
}
where $\eps,\eta_{1,2}$ are parameters, $\eta_d=\eta_2-\eta_1$, 
$W$ is a double--well--potential, 
typically 
containing more parameters, with $W'(-1)=W'(0)=W'(u_+)=0$, for some 
$u_+>0$, and ${\cal G}$ is an operator ensuring mass-conservation, 
e.g., $\CG f=f-\frac{1}{|\Om|}\int_\Om f(x)\dd x$. 
In suitable parameter regimes \reff{tfch} 
is extremely rich in pattern formation. The basic building blocks are 
straight (see Fig.\ref{fchf1})
and curved  (see Fig.\ref{fchf1c} and Fig.\ref{fchf1b}) ``channels'', 
i.e., bilayer 
interfaces between $u\equiv -1$ and some positive $u$, which 
show ``pearling'' and ``meander'' instabilities, leading to more 
complex patterns. 

Here we explain how to start exploring these with p2p2. 
Setting $v=\eps^2\Del u-W'(u)$, the stationary equation can be written 
as the two component system 
\begin{subequations}\label{sfchs}
\hual{
&-\eps^2 \Delta u+W'(u)+v=0,\\
&-\eps^2 \Delta v+W''(u)v-\eps\eta_1 v-\eps\eta_dW'(u)+\eps\ga=0, 
}
where $\ga$ is a Lagrange-multiplier for mass-conservation in \reff{tfch}. 
We take $\ga$ as an additional unknown,  and add the equation 
\huga{\label{sfchs3}
q(u):=\int_\Om u\dd x-m=0, 
}
\end{subequations}
where $m$ is a reference mass, also taken as a parameter.  
Thus, we now have 4 parameters $(\eta_1,\eta_2,\eps,m)$, one 
additional unknown $\gamma$, and one additional equation $n_q=1$. 
To implement $q$ from \reff{sfchs3}, and, 
strongly recommended, also $\pa_u q$, we set 
{\tt fuha.qf=@fchqf; fuha.qfder=@fchqjac;} and {\tt sw.qjac=1}. 
For $W$ we follow \cite[\S5]{DHPW12} and let 
$$
W(u)=W_p(u+1)+20(u-m_p+1)^{p+1}H(u-m_p+1),\text{ where }
W_p(u)=\frac 1 {p-2}(pu^2-2u^p)
$$
with $m_p=(p/2)^{1/(p-2)}$, 
and $H$ being the Heaviside function. In \cite[\S5]{DHPW12} numerical 
time integrations are presented with $p=3$, $\eps=0.1$, 
$\eta_2=2$, and $\eta_1=1$ (which leads to pearling), 
resp.~$\eta_1=2$ (which gives meandering). We aim at similar parameter 
regimes, but remark that we use somewhat larger $\eps$ to keep 
numerical costs low in our tutorial setting.

Getting good initial guesses for continuation is a delicate problem 
for \reff{sfchs}. Here we use guesses of the form 
\huga{\label{fchiguess}
\text{$u_{ig}(x)=-1+(a_1+a_2\sin(b_1 y)/\cosh(b_2x), \quad 
v_{ig}=-W'(u_{ig})$}} 
on rectangular domains $|\Om|$ with homogeneous Neumann BC for $u$ and $v$, 
and regular initial meshes of about 8000 points. 
Moreover, we can either let the software calculate the mass 
$m=\frac 1 {|\Om|}\int u(x)\dd x$ of 
the initial guess and use it as the constraint, or give a target $m$ externally. 

Figures \ref{fchf1}(a)-(c) show some first continuation of and bifurcation from 
a straight channel, for $\eps=0.35, \eta_2=2.5$, fixing $m=-0.9$, 
primary parameter $\eta_1$, and initial guess at $\eta_1=1$ 
of the form \reff{fchiguess} with 
$a_1=1.25, a_2=0, b_1=0, b_2=10$. 
The interfaces between $-1$ and $u_+$ and vice versa are rather 
sharp already for $\eps=0.35$, and only via adaptive mesh--refinement 
we get a ``straight channel'' 
solution from \reff{fchiguess}, on a grid of 
about 15.000 triangles. Then increasing 
$\eta_1$ we get a number of ``pearling instabilities'' 
and bifurcating pearling branches, which moreover show 
secondary bifurcations. Interestingly, while the 2nd and 4th 
primary pearling instabilities have roughly the same wavelengths 
as the first (see (b)), the 3rd has a much shorter wavelength (see (c)). 

As indicated above, slightly changing initial guesses may lead to 
failure of the initial Newton loop, or to convergence to quite 
different solutions. As an example we present in Fig.~\ref{fchf1}(d) 
the solution for $\eps=0.35, \eta_1=1, \eta_2=3$ 
from an initial guess of the form \reff{fchiguess} with 
$a_1=1.25, a_2=0.25, b_1=6, b_2=10$, and requiring $m=-0.98$. 
This directly yields a pearled straight channel. 

Getting a curved channel from (a curved version of) an initial guess 
like \reff{fchiguess} turns out to be difficult for small $\eps$. 
Figure \ref{fchf1c} shows an example obtained after some trial and error 
(where the Newton loop often gets stuck at residuals of about 
$10^{-4}$ or $10^{-5}$), after mesh--refinement to about 
25.000 triangles, but with still rather large $\eps=0.75$. 
After having found some curved solution (on the black branch), we can continue it to smaller $\eps$, which, however, needs more mesh refimement, and the bifurcation scenario does not change 
much compared to the case $\eps=0.75$ in Fig.~\ref{fchf1c}.

Finally, in {\tt fchcmds3.m} we give firstly
a curved channel obtained by first using 
\mlab's {\tt fsolve} on an initial guess, see Fig.~\ref{fchf1b} (a), 
and secondly a channel with a sharp bend as an initial guess 
that leads to spots, Fig.~\ref{fchf1b} (b). 

\brem\label{ntotrem}{\rm 
As the {\tt fCH} examples are rather slow, in the {\tt cmds} files 
we often set {\tt p.nc.ntot} to a rather 
small number, e.g., {\tt p.nc.ntot=20}, and we often switch off 
bifurcation detection and localization on bifurcating branches. 
When using these files as templates, this very likely needs to be reset.}
\eex\erem 


\subsection{Phase equation for traveling waves}\label{s:tw}
Systems \eqref{gform} with continuous symmetries require the selection
of a particular group element in order to allow for a continuation approach since
otherwise the linearization $\partial_u G$ has a kernel. 
This can often conveniently be
done by adding a suitable constraint, such as the norm constraint in 
\eqref{heq} for the scaling symmetry in the eigenvalue problem. 
 Note that although
the discretization error may eliminate the kernel, such that the continuation 
works, it is
essentially uncontrolled and correct resolution of parameter
dependencies requires a modification. 

A traveling wave $u$ on an infinite strip $\Omega=\R\times(-L,L)$
possesses a translation symmetry in $x$-direction with kernel of
$\partial_u G(u)$ generated by the spatial derivative $\partial_x
u$. The selection of a fixed translate is naturally done by (i) adding
a comoving frame term $s \partial_x u$ to \eqref{gform}, that is, by
modifying the tensor $b$ with an additional parameter $s$, and (ii)
adding an auxiliary equation that constrains the continuation path to
be orthogonal to the group orbit: $\langle \partial_\lambda
u, \partial_x u\rangle_2=0$. Numerically, we discretize the derivative
in the continuation direction $\partial_\lambda u \approx (u_{\rm
  old}-u)/(\lambda_{\rm old} - \lambda)$, where the subscript `old'
refers to the previous continuation step (saved in \texttt{p.u} in
\pdepb). Since the division by $\lambda_{\rm old} - \lambda$ is
redundant, we obtain the auxiliary equation 
$$\langle \partial_x u,
u_{\rm old}-u\rangle_2=0.$$ 
For the benefit of a simpler derivative of the resulting $G$ with respect to $u$ one may also use $\langle \partial_x u_{\rm old}, u_{\rm old}-u\rangle_2=0$ since these are equivalent in the continuum limit. In the case of periodic boundary conditions in the $x$-direction we have $\langle \partial_x u, u\rangle=0$ such that the conditions simplify to 
$$\text{$\langle \partial_x u, u_{\rm old}\rangle=0$ or 
$\langle \partial_x u_{\rm old}, u\rangle=0$,}
$$ 
respectively. 
See, e.g., demo \texttt{schnaktravel} for an implementation in \pdepb.

On the other hand, with the separated boundary conditions in \mlab's 
\ptool\ the translation symmetry never appears, but the constraining
procedure still allows to model the real line. For illustration,
consider fronts, which are spatially heteroclinic connections to
homogeneous steady states. In this case a sufficiently long  
$\Omega$ in the $x$-direction and homogeneous Neumann \bcs \ at the end
sides is (generically) a small perturbation to the front profile, also
for a traveling front with $s$ chosen as the wave speed. Continuation
in a parameter $\lambda$ of $G$ is now typically possible, which means
\emph{fixed} speed $s$. But on the infinite strip the speed $s$ will
typically depend on $\lambda$. This is resolved by a constraint as
described above, which (implicitly) couples $s$ and $\lambda$, so that
continuation of the extended system in $\lambda$ with additional
unknown $s$ calculates their interdependence. In \S\ref{pbcsec} 
we show how periodic domains can be implemented in \pdepb\ such that
translation symmetry can be realized (e.g. for periodic or localized traveling waves) and a priori requires the
constraint.

\paragraph{Example: Traveling fronts in an Allen-Cahn model ({\tt acfront}).}

For illustration of the simplest setting, consider the Allen-Cahn equation $\partial_t u =
\Delta u + \lambda u(1-u)(\mu+u)$, whose traveling waves with speed $s$ in $x$-direction solve the elliptic equation
\[
-\Delta u - \lambda u(1-u)(\mu+u) - s \partial_x u = 0.
\]
The (explicitly known) quasi 1D traveling front solutions are near a $y$-independent heteroclinic connection from $-\mu$ to $1$. With domain $\Omega$ set to a finite rectangle with homogeneous
Neumann \bcs~ we detect near-front solutions as follows. For $\mu=1$ the nonlinearity is symmetric and waves
stationary. Increasing $\lambda$ from $1$ yields a pitchfork
bifurcation to half cosine-modes that approach a heteroclinic
connection between $\pm 1$ as $\lambda$ increases further. Next we add
the constraint and the parameter $s$ and perform a continuation in the
symmetry breaking parameter $\mu$. The bifurcation diagrams (see 
{\tt acfront\_cmds})
match the explicitly known dependence of $s$ on $\mu$, e.g.\ \cite{grind}.

\subsection{Periodic boundary conditions for rectangular 
domains}\label{pbcsec}

For axis-aligned rectangular domains \pdepb~ can identify opposite sides with equal
grid arrangements\footnote{\pdepb~only checks the coordinate value in
  the periodic direction and assumes equal number of points; the
  transverse direction can have other shapes.} in order to generate
cylindrical or toroidal geometry. The initial setup requires
homogeneous Neumann boundary conditions on the sides that are to be
identified, and the grid requirement is most easily realized with a
mesh from \texttt{poimesh}. The boundary conditions on the remaining
boundary can be arbitrary. For all calculations the effective mesh is
reduced by removing the points from one of the identified sides of the
rectangle such that the solution vector \texttt{p.u} is smaller than
on the initial mesh. However, the full mesh and the Neumann BC are used
for assembling the FEM discretized PDE and for plotting purposes.

The transformation of a vector from the reduced to the full mesh goes
by the matrix \texttt{p.mat.fill}, which simply extends a vector by
generating copies of entries on the periodic boundary. For instance,
\texttt{p.mat.fill*p.u(1:p.nu)} gives the extended solution vector.

The switch to periodic domains in \texttt{p2p2} can be conveniently
done by calling the routine \texttt{rec2per} with additional argument
to determine the type of periodic domain: 

\medskip
1: top=bottom side, 2: left=right side, 3: torus (this setting is stored in \texttt{p.sw.bcper}). 

\medskip
The convenience function \texttt{rec2perf} in
addition loads a point from a Neumann \bcs~ solution from a file for the purpose of continuing from this solution with periodic \bcs.

These routines essentially generate the matrix \texttt{p.mat.fill},
which is set to $1$ in the standard (non-periodic) setting. In addition
these routines modify a given solution from the Neumann domain to the
periodic domain by removing the redundant entries from \texttt{p.u}
with the matrix \texttt{p.mat.drop} via \texttt{p.mat.drop*p.u(1:p.nu)}; also
the degree of freedom parameter \texttt{p.nu} is set to the corresponding
(smaller) value.

Next, we explain the details of the transformation of the system matrices from Neumann to periodic \bcs\ and provide an
example in \S\ref{s:schnatw}. 

\brem\label{rem:saving2}{\rm 
Remark~\ref{frefrem} also applies to these matrices {\tt p.mat.XX} so that they are not saved to disk. In order to avoid miscalculations upon reloading, for periodic geometries \texttt{p.mat.fill} and \texttt{p.mat.drop} are saved as empty arrays {\tt []}  (for non-periodic domains these are saved as $1$); an attempt to run a problem from a periodic domain without resetting the matrices will then produce an error message. The matrices are automatically regenerated when loading a point with {\tt loadp}, which uses the geometry type stored in {\tt p.sw.bcper}.}
\eex\erem 

\brem{\rm 
Grid adaption for periodic domains is so far implemented only in a simple 
ad hoc way: to ensure that identified boundaries match after grid adaption, we 
remove triangles from a refinement--list generated by {\tt pdeadworst}. This is controlled 
by {\tt p.nc.bddistx, p.nc.bddisty}, see {\tt rmbdtri.m} for details. 
Thus, for periodic domains, mesh adaption is useful only as long as there are no large gradients near the identified domain boundaries. }
\eex\erem

\subsubsection{Transforming the FEM problem from Neumann to periodic \bcs.}

For illustration, consider the simple situation of a one-dimensional
chain with 4 elements. In that case node 1 and 4 are identified to
generate a ring and we have
\[
\texttt{p.mat.fill} = \begin{pmatrix}
1 & 0 & 0\\
0 & 1 & 0\\
0 & 0 & 1\\
1 & 0 & 0
\end{pmatrix}, \quad
\texttt{p.mat.fill'} = \begin{pmatrix}
1 & 0 & 0 & 1\\
0 & 1 & 0 & 0\\
0 & 0 & 1 & 0
\end{pmatrix}, \quad
\texttt{p.mat.drop} = \begin{pmatrix}
1 & 0 & 0 & 0\\
0 & 1 & 0 & 0\\
0 & 0 & 1 & 0
\end{pmatrix}.
\]
Hence, \texttt{p.mat.fill} writes a copy of entry 1 into slot 4 while
\texttt{p.mat.fill'} adds entries 1 and 4 into slot 1, and
\texttt{p.mat.drop} simply removes the last entry. Note that in the
actual matrix construction it is not assumed that the points that are
to be identified appear within \texttt{p.u} in any specific ordering.

\medskip Next, observe that the piecewise linear `hat' basis function
at a node is the sum of its triangular parts over the neighboring
triangles. At a boundary node the basis function for homogeneous
Neumann conditions simply does not have an additive contribution from
outside the grid. Therefore, it can be extended to a full basis
function of an interior node in the periodic domain by adding the
corresponding contributions from periodically identified nodes. Denote the basis functions for the Neumann problem  by $(\phi_j)_{j=1}^{n_p}$ and those for the periodic problem by $(\psi_j)_{j=1}^{n_p^\text{per}}$ with $n_p^\text{per} < n_p$. For the above 1D example $n_p=4, n^\text{per}_p=3, \psi_1 = \phi_1+\phi_4, \psi_2=\phi_2$, and $\psi_3=\phi_3$. Considering for simplicity the mass matrix (denoted by $M$in the Neumann case and $M^\text{per}$ in the periodic case), we have $M_{ij}=\langle \phi_i,\phi_j\rangle_2$  and $M^\text{per}_{ij}=\langle \psi_i,\psi_j\rangle_2$ such that
\begin{align*}
&M^\text{per}_{ij}=M_{ij} \quad \text{for} \ 2\leq i,j\leq 3\\
&M^\text{per}_{1j}=M_{1j}+M_{4j}, \ M^\text{per}_{j1}=M_{j1}+M_{j4},  \quad \text{for} \ 2\leq j\leq 3\\
&M^\text{per}_{11}=M_{11}+M_{14} +M_{41} + M_{44}. 
\end{align*}
The modification of the Neumann stiffness matrix $K$ to the periodic stiffness matrix $K^\text{per}$ is completely analogous. The right hand side $F$ satisfies
$$F^\text{per}_j = F_j \quad \text{for} \ 2\leq j\leq 3, \ F^\text{per}_1 = F_1+F_4.$$
These transformations are efficiently performed via
\begin{center}\begin{tabular}{rl}
\texttt{K}$_\per$ &= \texttt{p.mat.fill' * K * p.mat.fill}\\
\texttt{M}$_\per$ &= \texttt{p.mat.fill' * M * p.mat.fill}\\
\texttt{F}$_\per$ &= \texttt{p.mat.fill' * F}.
\end{tabular}\end{center}
In practice we assemble $K,M$ and $F$ via the \mlab~\ptool\ and apply the above transformation.

In order to account for the dependence of $c,b,a, f$ in \eqref{gform} on \texttt{p.u}, the vector \texttt{p.mat.fill*p.u} is fed into the Neumann BC assembling routines. Recall that the matrices {\tt p.mat.fill} and {\tt p.mat.drop} are generated automatically by the routine {\tt rec2per} according to the value of {\tt p.sw.bcper}.

\subsubsection{Cylinder geometry: Quasi-1D traveling waves 
({\tt schnacktravel})}\label{s:schnatw}
For the simplest illustration of cylinder geometries 
we consider a quasi 1D setting in the Schnakenberg model \eqref{mod1} with
$\sigma=0$, $\Omega=(-0.1,0.1)\times(-L,L)$ and full Neumann boundary conditions at
first. Based on expectations of the general structure of the existence
region for wavetrains in such systems  \cite{JRsemi, GSBusse}, we follow a certain continuation path in
$(\lambda,L)$ that leads to a travelling wave bifurcation. 
First we choose the domain length $L=L_c$ compatible with the Turing instability mentioned in \S\ref{ssec}: this occurs at $\lambda_c = \sqrt{d}\sqrt{3-2\sqrt{2}}\approx 3.2085$ with spatial period $L_c=2\pi/k_c$, $k_c= \sqrt{\sqrt{2}-1}$.

\begin{figure}
\begin{center}
\begin{tabular}{llll}
(a)&(b)&(c)&(d)\\
\ig[width=0.22\textwidth]{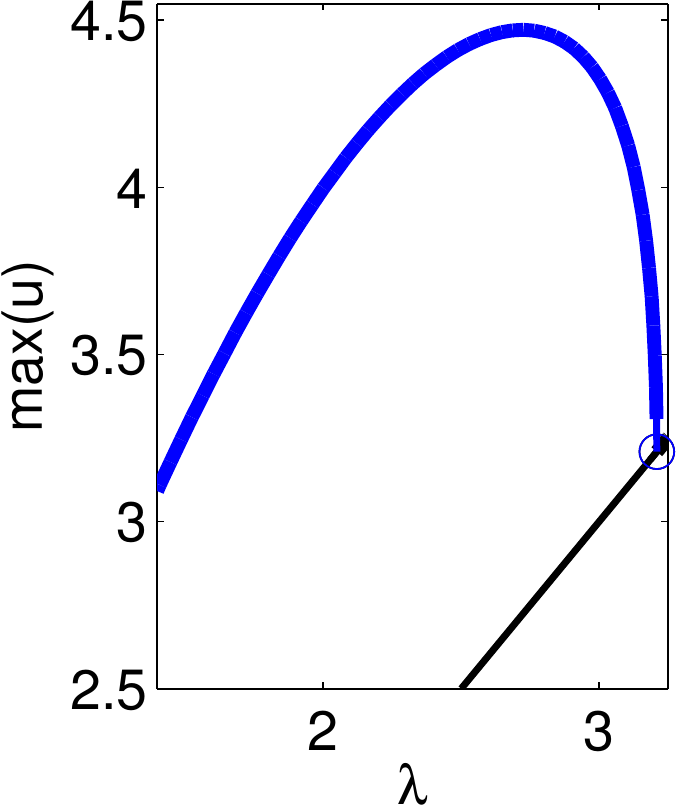} &
\ig[width=0.22\textwidth]{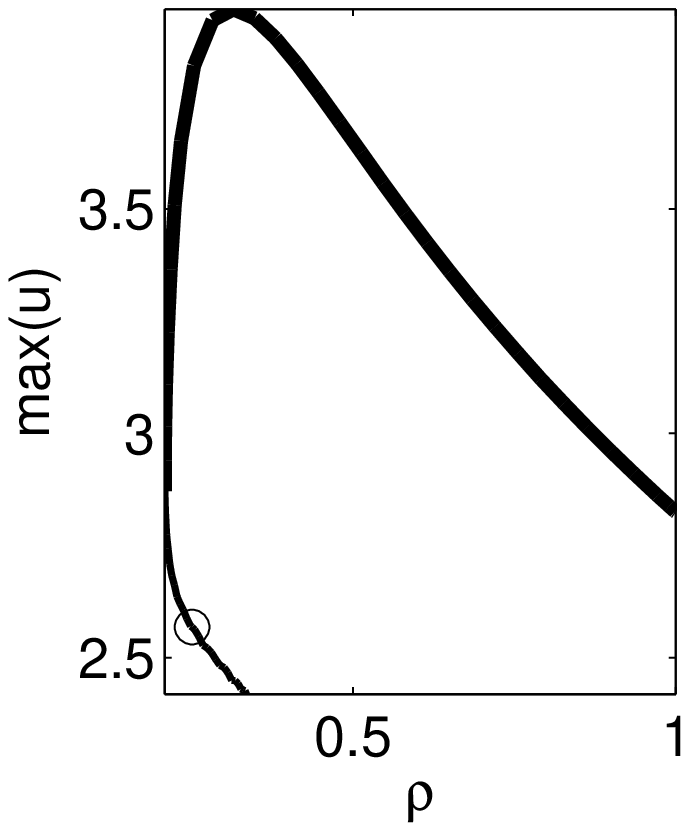}&
\ig[width=0.22\textwidth]{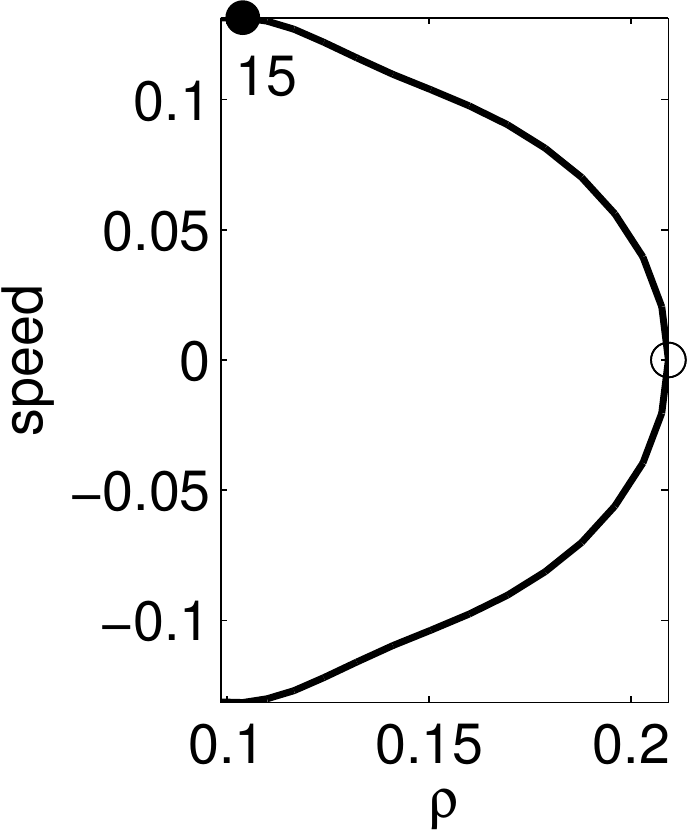} &
\ig[width=0.22\textwidth]{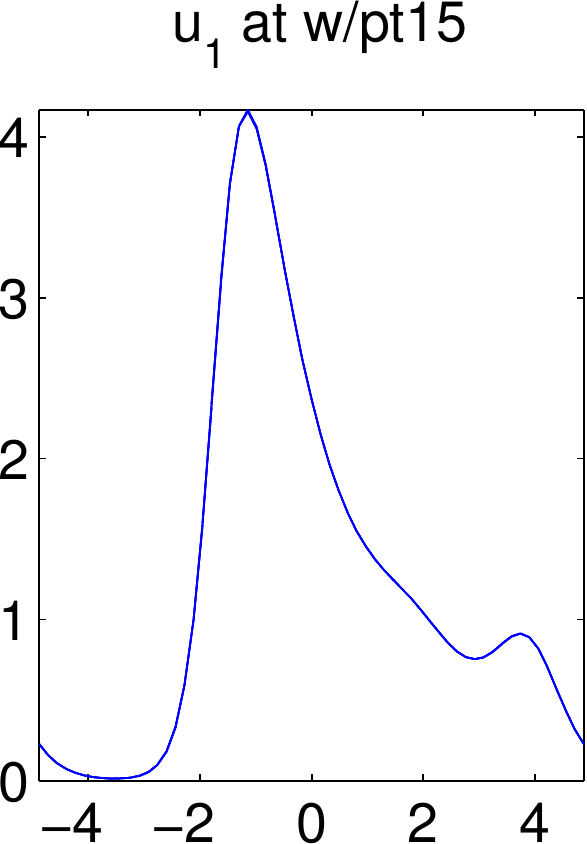}
\end{tabular}\end{center}
\caption{(a) Bifurcation of stripes under Neumann \bcs~ through a 
Turing instability. (b) Continuation from endpoint of (a) in 
the domain length parameter $\rho=(L_c/L)^2$. Here the circle marks a bifurcation under Neumann \bcs. The traveling wave bifurcation under periodic \bcs~ is 
closer to the fold point. (c) Bifurcation diagram in cylindrical geometry with parameters $\rho$ and speed $s$. (d) profile for the solution at 
point 15 marked in (c).
\label{f:schnaksteady}}
\end{figure}

For decreasing $\lambda$ from, say $\lambda=3.5$, the instability appears as a pitchfork bifurcation under Neumann BC and we continue the resulting branch until $\lambda=1.3$. See Figure~\ref{f:schnaksteady}(a). This value is somewhat arbitrarily chosen from the expectation that for sufficiently small $\lambda$, increasing the domain length yields a fold and traveling wave bifurcation for the associated wavetrains.  To change the domain length,  we multiply both diffusion constants by a factor $\rho$, i.e. $L=\pi/(\sqrt{\rho} k_c)$. Now we switch the continuation parameter to $\rho$ and continuation starting from the endpoint of the previous continuation indeed leads to a fold. See Figure~\ref{f:schnaksteady}(b).  Note that
the bifurcation point marked with a circle does not give traveling waves on $\R$ as it stems from Neumann \bcs. This is because there is no extension of the Neumann solution from the bounded domain onto $\R$.

After these preparatory steps we change the geometry to a cylinder
such that the boundaries at $-L$ and $L$ are identified. As described
in \S\ref{s:tw}, we add a phase equation to eliminate the zero
eigenvalue from translation symmetry and add the traveling wave speed as a second parameter to the
active continuation parameters. Finally, we load the endpoint from the
previous calculation. These steps are conveniently done in the
{\tt schnaktravel} demo by the following commands:

{\small\begin{verbatim}
r=rec2perf('q2','pt25','r',1);  % load 'q2/pt25', set directory name 'r', geometry type 1
r.nc.ilam=[3;2];                % set parameters with indices 3 and 2 as active
r.nc.nq=1; r.fuha.qf=@schnakqf; % set number auf aux. eqn. and function handle
p.sw.qjac=1; p.fuha.qfder=@schnak_qfder;         % analytical jac for aux. eqn.
\end{verbatim}
}
Here the new primary parameter has index 3, which corresponds to the traveling wave speed $s$ in the comoving frame term $s \partial_y(u,v)^T$ added to \eqref{gform}; see \S\ref{s:tw}. 
\medskip Now we perform a continuation from the endpoint of the Neumann \bcs\ computations back for decreasing $\rho$. The stationary solution branch is the same, but the location of the bifurcation point changed to a value much closer to the fold point. Branch switching in both
directions and continuation yields the branches and profile plotted in Figure~\ref{f:schnaksteady}(c),(d). 

\subsubsection{2D traveling waves in a cylinder ({\tt twofluid})}\label{S:2fluid}

A traveling wave problem that cannot be reduced to a one dimensional problem occurs in the `two-fluid' tokamak plasma model from \cite{tf}. The profile 
satisfies the elliptic problem 
\def\DT{\delta}
\begin{equation}\label{e:tf}
\begin{aligned}
0 &= - \nu \Delta \uP - (\nabla V)^{\perp} \cdot \nabla \uP -(\DT+s) \partial_{x_2} \uP + \partial_{x_2}V/L_1 ,\\
0 &= - \nu \Delta \uM - (\nabla V)^{\perp} \cdot \nabla \uM  -s \partial_{x_2} \uM - \partial_{x_2}V/L_1 ,\\
0 &= -\Delta V -\uP - \uM 
\end{aligned}
\end{equation}
posed on $\Omega=(-1,1)^2$ subject to periodic \bcs~ in the second variable $x_2$ and homogeneous Dirichlet \bcs\ in $x_1$, with parameters $\nu,\del,L_1>0$. 
For definiteness we fix $L_1=2$ and $\nu=9\cdot 10^{-4}$, and take 
$\del$ as primary parameter. 

\begin{figure}
\begin{center}
\begin{tabular}{llll}
(a)&(b) &(c) &(d)\\
\ig[width = 0.22\textwidth]{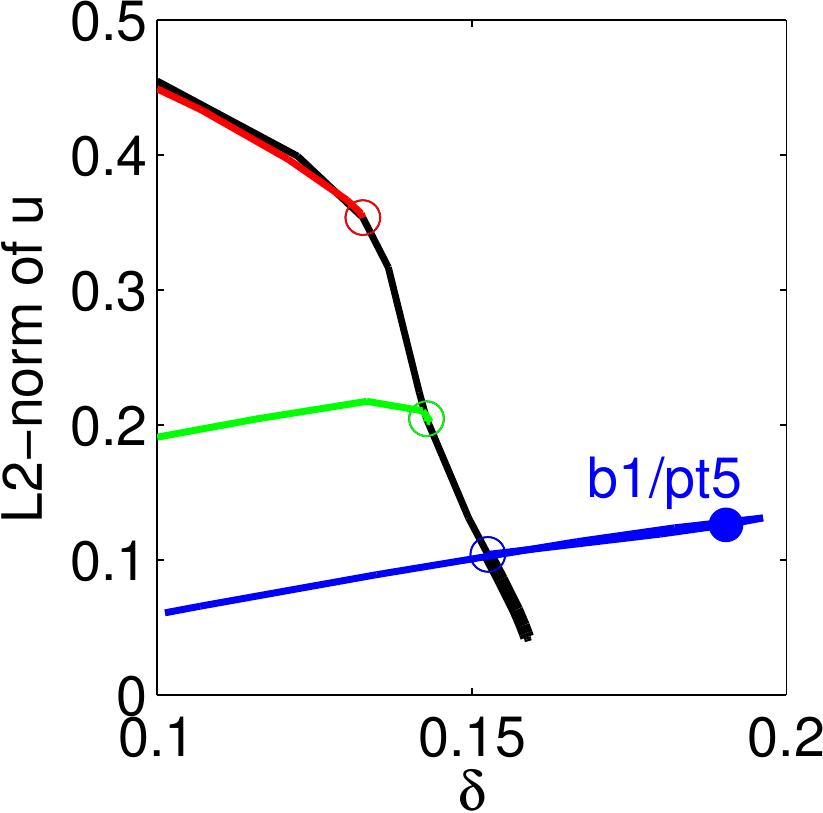} & 
\ig[width = 0.21\textwidth]{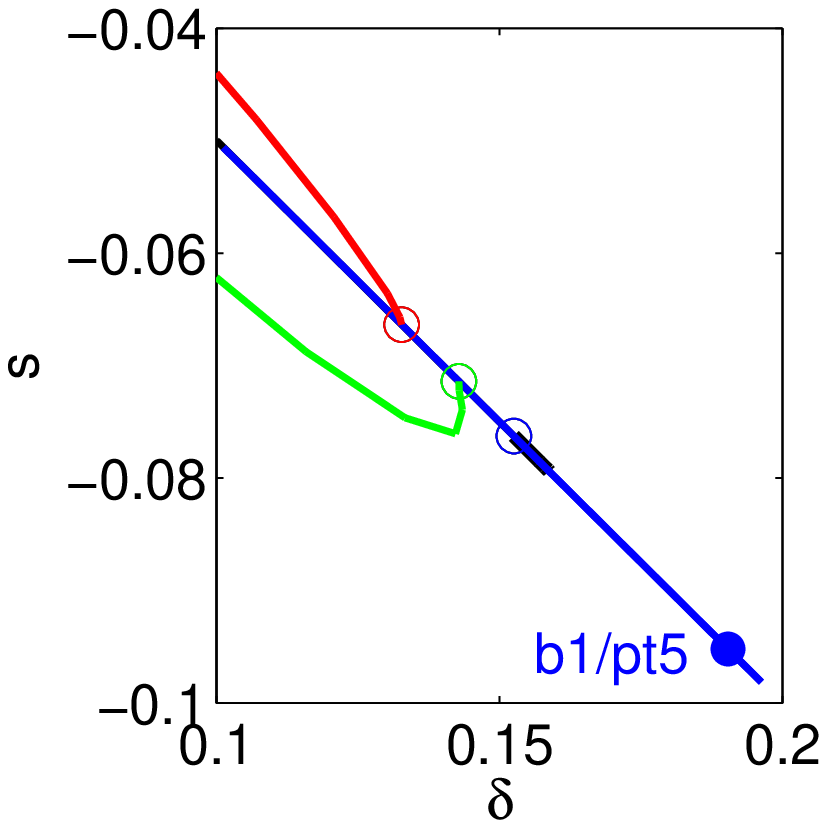} &
\ig[width = 0.25\textwidth]{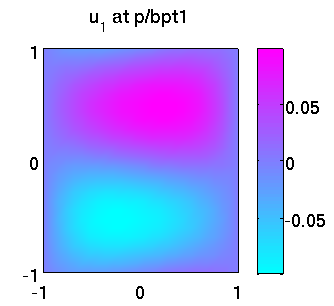}&
\ig[width = 0.25\textwidth]{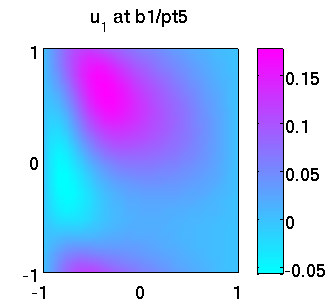}
\end{tabular}
\end{center}
\caption{(a),(b) Bifurcation diagrams. (c),(d) Contour plots of $u_1$-component 
for solutions.\label{f:tf1}}
\end{figure}

As proven in \cite{tf}, in the parabolic problem the trivial state $u_1=u_2=V=0$ for $s=0$ undergoes a generic supercritical Hopf bifurcation for decreasing $\DT$ at $\DT=\DT_2\approx 0.16$ and this generates periodic traveling waves. However, using the knowledge of $s$ at onset, the bifurcation has a double zero eigenvalue and hence cannot be detected by the standard {\tt bifdetec}. Instead, we set $s$ at the onset speed and $\DT$ just below the bifurcation value, and perturb the trivial solution by a rough eigenfunction. 
Using time evolution with {\tt tint}, the trajectory indeed approaches the stable traveling  wave solution. A Newton-loop for the system augmented by the phase equation from \S\ref{s:tw} then yields an initial solution to \eqref{e:tf}. In Figure~\ref{f:tf1} we plot some resulting branches and solutions.

Remarkably, we find solutions on the blue branch for values of $\DT$ \emph{larger} than the bifurcation point. It thus coexists with the \emph{locally} stable trivial state, which is proven to be \emph{globally} stable  for $\DT>5\DT_2$ in \cite{tf}. This numerical result thus shows that failure to prove global stability of the trivial state for $\DT\in(\DT_2,5\DT]$  is not a technicality: there is a branch of nontrivial solutions preventing global stability for some range above $\DT_2$. 

\subsubsection{Torus geometry: nonlinear Bloch waves ({\tt nlb})}\label{S:NLB}
As a second example for periodic boundary conditions we consider the
time harmonic Gross-Pitaevskii equation
\begin{equation}\label{E:SGP}
  -\Delta \phi -\omega \phi +V(x) \phi +\sigma |\phi|^2 \phi =0, \quad x\in \R^2
\end{equation}
with the periodic potential $V(x+2\pi e_m) = V(x)$ for all $x\in \R^2$
and $m=1,2$, where $e_m$ is the $m$-th Euclidean unit vector in
$\R^2$, and $\sig=\pm 1$. Equation \eqref{E:SGP} describes, e.g., time harmonic
electromagnetic fields in nonlinear photonic crystals or Bose-Einstein
condensates loaded on optical lattices. It has quasi--periodic
solutions bifurcating from the trivial solution at spectral points
$\omega_* \in \text{spec}(-\Delta +V)$, see \cite{DPS09,DU14}. Here we consider
the particular case where
$\omega_*=\omega_{n_*}(k_*)$ for some fixed $k_*\in (-1/2,1/2]^2$, 
with $\omega_{n_*}$ the $n_*$-th band function in the band
structure $(\omega_n(k))_{n\in \N}$, $k\in (-1/2,1/2]^2$. Thus we seek a
quasi-periodic nonlinear Bloch-wave $\phi$ of \eqref{E:SGP} with the
quasi-periodicity vector $k_*$, i.e., 
$\phi(x) =e^{\ri k_*\cdot x} \eta(x), \eta(x{+}2\pi e_m)=\eta(x)\text{ for all } x{\in}\R^2, m=1,2.$
As shown in \cite{DU14}, for $\omega=\omega_*+
\text{sign}(\sigma)\varepsilon^2$ with $\varepsilon >0$ small enough
such nonlinear Bloch waves exist and have the asymptotics
\huga{\label{nlbas}
\text{$\phi(x) \sim \varepsilon
\left(|\sigma|\|p_{n_*}(\cdot,k_*)\|_{L^4((-\pi,\pi)^2)}^4\right)^{-1/2}
p_{n_*}(x,k_*)e^{\ri k_*\cdot x}$  for $\eps \to 0.$}}
Inserting $\phi(x) =e^{\ri k_*\cdot x} \eta(x)$ in \eqref{E:SGP} and
using real variables $u_1,u_2$, i.e., $\eta=u_1+\ri u_2$, we
get
\huga{\label{SGPb}
0=G(u_1,u_2):=-\begin{pmatrix} \Delta u_1 \\ \Delta u_2\end{pmatrix}+2 \begin{pmatrix} k_*\cdot \nabla u_2\\ -k_*\cdot \nabla u_1\end{pmatrix} + (|k_*|^2-\omega +V(x)) \begin{pmatrix} u_1 \\ u_2 \end{pmatrix} +\sigma (u_1^2+u_2^2) \begin{pmatrix} u_1 \\ u_2 \end{pmatrix}, 
}
on the torus $\mathbb{T}^2=\R^2/(2\pi\Z^2)$, with linerization 
\hugast{
Lu=-\begin{pmatrix} \Delta u_1 \\ \Delta u_2\end{pmatrix}+2 \begin{pmatrix} k_*\cdot \nabla u_2\\ -k_*\cdot \nabla u_1\end{pmatrix} + (|k_*|^2-\omega +V(x)) \begin{pmatrix} u_1 \\ u_2 \end{pmatrix}
} 
around $0$. We aim to find 
bifurcations from the branch $u{=}0$, with primary 
parameter $\om$. 

\reff{SGPb} has the continuous symmetry $u\mapsto 
\ssmatrix{\cos\al&-\sin\al\\\sin\al&\cos\al}u$ from the phase 
invariance $\phi\mapsto\er^{\ri\al}\phi$ of \reff{E:SGP}. In particular, 
each eigenvalue 
$\mu$ of $L$ is double: if $Lu=\mu u$, then also, e.g., $Lv=\mu v$ with 
$v=(-u_1,u_2)$ which is linearly independent of $u$. 
To deal with this phase invariance we proceed as follows. 
First we modify {\tt findbif.m} to search for $\om$ values 
where {\em two} eigenvalues go through zero; see {\tt findbifm.m}, 
which can easily be generalized to $m$ eigenvalues going through zero. 
Next we modify {\tt swibra.m} to a version {\tt swibram.m}. 
Here, the user first has to choose which of the two 
zero eigenvalues to use for bifurcation, and second, has to set a 
``phase-fix-factor'' ${\tt pffac}=\gamma\ne 0$. 

For $\gamma=0$, branch--switching proceeds as usual: {\tt swibram} produces a new tangent $\tau$, 
from which we may continue a bifurcating branch using {\tt cont}, 
with one caveat. During continuation, the phase $\al$ of the solution, 
for instance defined as $\al=\arctan\frac{u_2(x_*)}{u_1(x_*)}$ at 
one point $x_*$ in the domain, 
may change in an uncontrolled way, strongly dependent, for instance, 
on the step--length. \footnote{That continuation works at all, 
despite the zero eigenvalue from the phase--invariance of \reff{SGPb}, 
is due to the Fredholm alternative: in the Newton loops, the RHS 
is perpendicular to the kernel. See also \cite[\S5.1]{p2p}.} 

For $\ga\ne 0$, the software removes this undesired phase--wandering by fixing $u_2(x_*)=0$ at a point $x_*$ where $|u_1(x_*)|$ (i.e., the absolute value of the real part) 
is maximal. 
If $\gamma <0$, this is achieved by the method (a), where $u_2(x_*)$ is set to $0$ and $u_2(x_*)$ as well as the 
equation for $u_2(x_*)$ are dropped from the discretized system. For this we use modifications of {\tt p.mat.drop} and {\tt p.mat.fill}, 
see {\tt dropp.m}. 

If $\gamma >0$ is chosen, the following method (b) is used. 
Assume that $x_*$ is point ${\tt p.pfn}=n_*$ in the 
discretization. We then overload {\tt pderesi.m} locally to 
replace the $n_u/2+n_*$ entry  {\tt r(p.nu/2+p.pfn)} 
in $G(u)$ by $\ga u(n_u/2+n_*)$, i.e., 
{\tt r(p.nu/2+p.pfn)=p.pffac*u(p.nu/2+p.pfn);} 
where {\tt p.pffac}=$\gamma$. Accordingly we also overload {\tt getGupde.m} locally and add:  \\
{\tt Gu(p.nu/2+p.pfn,:)=zeros(1,p.nu); Gu(p.nu/2+p.pfn,p.nu/2+p.pfn)=p.pffac;}

In {\tt swibram.m} we thus use method (a) if ${\tt p.pffac}<0$ (with the precise value 
irrelevant) and method (b) if ${\tt p.pffac}>0$ (with for instance $\ga=10^3$). 
Both methods yield indistinguishable results for all our tests. 

\begin{figure}
  \begin{center}
\ig[height=40mm,width=35mm]{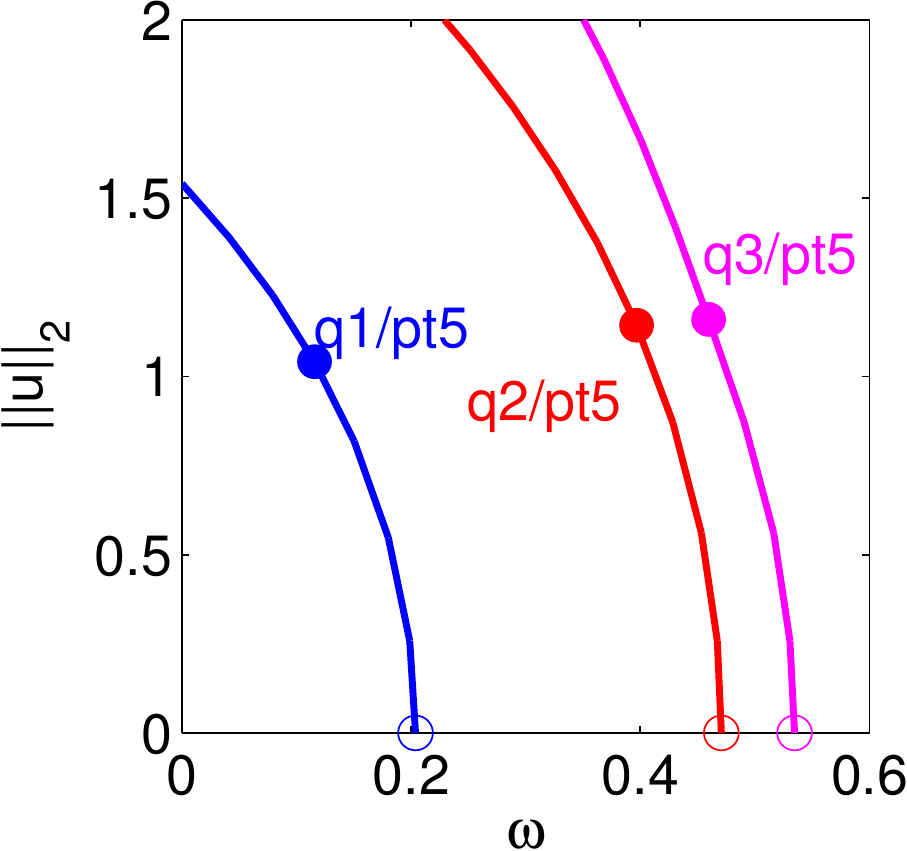}
\ig[height=40mm]{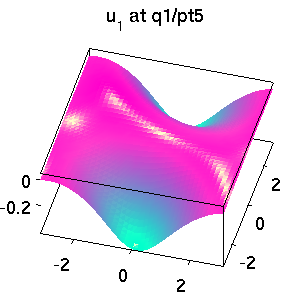}
\ig[height=40mm]{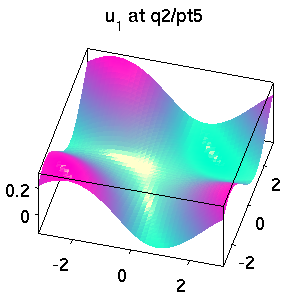}
\ig[height=40mm]{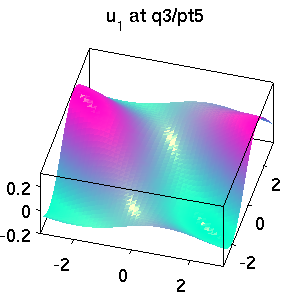}
\caption{{\small Nonlinear Bloch wave example for
      $V(x)=e^{-x_1^2}\cos(x_2)$ and $k_*=(1/2, 1/2)$. }\label{fnlb1}}
  \end{center}
\end{figure}

For a numerical example we choose $\sigma =-1$, the potential
$V(x)=e^{-x_1^2}\cos(x_2), \ x\in (- \pi,\pi]^2$ and 
$k_*=(1/2, 1/2)$.  Figure~\ref{fnlb1} shows the first three bifurcating 
branches, and real parts of selected profiles.  For all three 
branches we obtain excellent agreement with the asymptotics \reff{nlbas}. 
The mesh was generated by \texttt{poimesh}. 
Though time is not crucial here, for illustration 
we precompute ${\tt pot}:=V$ resp.~its interpolation 
{\tt poti} to the triangle centers (needed in mesh-refinement)  
and put these fields into 
{\tt p.mat}. As {\tt p.mat} is {\em not} saved to disk, we then also 
need to overload {\tt loadp.m} locally, and also recompute 
{\tt p.mat.pot} and {\tt p.mat.poti} after mesh--refinement, 
see {\tt nlbpmm.m}. See also the end of {\tt cmds.m} for an example 
of mesh-refinement following Remark \ref{rem:saving2}. 

\subsection{Other examples from \cite{p2p}}
The root demo directory {\tt p2p2/demos} also contains transfers of the 
examples from \cite{p2p} to the new setup, e.g: 
\bci 
\item {\tt acgc}: the (cubi--quintic) Allen--Cahn with Dirichlet BC and 
a global coupling, 
i.e., $$-0.1\Delta u-u-u^3+u^5-\lam\spr{u}=0,$$ where 
$\spr{u}=\int_\Om u\dd x$. This uses some modifications 
of the linear system solvers. 
\item {\tt bratu}: a scalar elliptic equation on the unit square with zero flux \bcs
\[
-\Delta u-f(u,\lam)=0, \quad f(u,\lam)=-10(u-\lam \er^{u}),
\]
for which a number of results can  be obtained analytically. The updated demo contains a branch point continuation, see  \S\ref{s:acfold}.

\item {\tt chemtax}: a quasi--linear non--diagonal reaction--diffusion system 
from chemotaxis in the form 
\huga{\label{chem1} 
G(u,\lam):=-\bpm D\Delta u_1-\lam \nabla\cdot(u_1\nabla u_2)
\\\Delta u_2\epm 
-\bpm r u_1(1-u_1)\\ \frac{u_1}{1+u_1}-u_2\epm=0. 
}
\item {\tt rbconv}: Rayleigh-B\'enard convection in the Boussinesq approximation 
streamfunction form
\begin{align}\label{e:rbconv}
  -\Delta\psi  + \omega&= 0, \nonumber\\
  -\sigma\Delta \omega -\sigma R \partial_x\theta + \partial_x\psi\partial_z\omega - \partial_z\psi\partial_x\omega &=0, \\
  -\Delta \theta - \partial_x\psi + \partial_x\psi\partial_z\theta
  - \partial_z\psi\partial_x\theta &=0,\nonumber
\end{align}
and with various boundary conditions. (The implementation given here detects both branches of the stress-free \bcs~by continuation -- in \cite{p2p} we used {\tt tint} to detect the second branch.)

\item {\tt gpsol}: 
time--harmonic Gross--Pitaevskii equations in a rotating frame, 
leading to real systems of the form 
\begin{subequations}\label{nls3}
\hual{
&-\Delta u+(r^2-\mu)u-|U|^2 u-\om(x\pa_y v-y\pa_x v)=0, \\
&-\Delta v+(r^2-\mu)v-|U|^2 v-\om(y\pa_x u-x\pa_y u)=0, 
}
\end{subequations} 
where $|U|^2=u^2+v^2$, and generalizations to more components. 
This is similar to \reff{SGPb}, but here we use potentials $V(x,y)=x^2+y^2$, 
and search for and continue soliton solutions. 
\item {\tt vkplate}: the Von K\'arm\'an equations for 
the buckling of elastic plates 
\huga{
\begin{split}
-\Delta^2 v-\lam\pa_x^2 v+[v,w]&=0, \quad 
-\Delta^2 w-\frac 1 2 [v,v]=0, 
\end{split}
\label{vk1} 
}
where $[v,w]:=v_{xx}w_{yy}-2v_{xy}w_{xy}+v_{yy}w_{xx}$, with various boundary 
conditions. After some transformations this yields a 
10-equations-system of the form \reff{gform}. 
\eci 
 In general, the transfer is rather straightforward, and here we only give 
some details on the new implementation of {\tt vkplate}
which is our most complicated example with respect to coding.  Mainly 
we want to illustrate how to gain additional flexibility in the 
{\tt sfem=1} setup by modifying 
{\tt setfemops.m}. This is needed if a parameter 
genuinely enters, for instance, the stiffness matrix $K$, as does $\lam$ in 
{\tt vkplate}. Similar ideas are also used in, e.g., {\tt gpsol} and {\tt rbconv}. See also Remark~\ref{r:cust}.

\subsubsection{Semilinear setting for Von K\'arm\'an equations 
({\tt vkplate})}\label{s:vkplate}
In \cite{p2p} we rewrite \reff{vk1} as a system of 10 equations, where the 
(bifurcation) parameter $\lam$ enters the stiffness matrix. In order to 
treat this in the {\tt sfem=1} setting, we essentially 
put the tensors {\tt a} and {\tt c} into {\tt p.eqn.a, p.eqn.c}, with one 
little trick: the $\lam$ dependent part is not put into {\tt p.eqn.c}, 
but (with $\lam=1$) into an extra field {\tt p.eqn.c2}. We then locally 
modify {\tt setfemops.m} to also assemble the associated stiffness matrix 
{\tt p.mat.K2}, and set up {\tt vksG} as, essentially, \\
{\tt r=(p.mat.K+lam*p.mat.K2)*u(1:p.nu)-p.mat.M*f;} 
The encoding of $\pa_u G(u)$ then also requires a little care, see the 
listing below. 

{\small\begin{verbatim} 
function p=setfemops(p) % modified for vkplate: generate additional K2 which 
% will be multiplied by lam in vksG and vksGjac 
upde=p.mat.fill*p.u(1:p.nu); neq=p.nc.neq; bc=p.fuha.bc(p,p.u); m=p.mesh; 
[~,p.mat.M,~,~,p.mat.bcG,~,~]=assempde(bc,m.p,m.e,m.t,0,1,zeros(neq,1),upde);
[p.mat.K,~]=assempde(bc,m.p,m.e,m.t,p.eqn.c,p.eqn.a, zeros(neq,1),upde);
[p.mat.K2,~]=assema(m.p,m.t,p.eqn.c2,0,zeros(neq,1));
end
\end{verbatim}}
{\small\begin{verbatim} 
function Gu=vksGjac(p,u)  % sfem=1 jacobian for von karman-plate
lam=u(p.nu+1); n=p.np; u5=u(4*n+1:5*n); u6=u(5*n+1:6*n); 
u7=u(6*n+1:7*n); u8=u(7*n+1:8*n); u9=u(8*n+1:9*n); u10=u(9*n+1:10*n); 
f2u5=spdiags(-u9,0,n,n); f2u6=spdiags(-u8,0,n,n); f2u7=spdiags(2*u10,0,n,n); 
f2u8=spdiags(-u6,0,n,n); f2u9=spdiags(-u5,0,n,n); f2u10=spdiags(2*u7,0,n,n);
f4u5=spdiags(u6,0,n,n); f4u6=spdiags(u5,0,n,n); f4u7=spdiags(-2*u7,0,n,n); 
zd=spdiags(zeros(n,1),0,n,n); % 0-diag for easy sorting of non-zeros into Fu
Fu=[sparse([],[],[],n,10*n,0); 
    [zd zd zd zd f2u5 f2u6 f2u7  f2u8  f2u9 f2u10]; 
    sparse([],[],[],n,10*n,0); 
    [zd zd zd zd f4u5 f4u6 f4u7   zd    zd   zd]; 
    sparse([],[],[],6*n,10*n,0)];  % set remainder of Fu to 0 
Gu=p.mat.K+lam*p.mat.K2-p.mat.M*Fu;
end
\end{verbatim}}

\section{Discussion and outlook}\label{todosec}
Compared to the version of \pdep\ documented in \cite{p2p}, \pdepb\ brings a 
number of 
\bcen[(i)]
\item extensions, e.g.: fold--and branchpoint continuation, 
general auxiliary equations, periodic boundary 
conditions,  interface to {\tt fsolve}, 
\item optimizations, e.g.: faster FEM in the {\tt sfem=1} setting, 
\item cleanups, reorganizations, improved user--friendliness, e.g.: 
substructures of {\tt p}, no single 
explicit parameter $\lam$ anymore, but easy switching between 
different parameters, improved plotting, 
\ecen 
and some bug-fixes (not documented in detail here). 
Moreover, besides the tutorial examples given here and in \cite{p2p}, 
\pdep\ and \pdepb\ have 
been applied to a number of genuine research problems, e.g., 
\cite{DU14, uwsnak14, uwfc14,tf}, and a number of further 
projects are in progress. Often, new projects require 
extensions of the software, and while (i) 
clears point 2 from the To-Do-List in \cite[\S6]{p2p}, the list 
as such has rather become longer. Currently, we are working on or 
planning the following extensions: 
\bcen
\item  Implement some more general (stationary) bifurcation handling, including 
branch--switching at multiple bifurcations; for the case of 
double eigenvalues due to phase--invariance this has been done 
in an ad--hoc way in \S\ref{S:NLB}. Also, Hopf bifurcations will be 
tackled. These points still roughly correspond to \cite[Point 1.~in \S6]{p2p}. 
\item Invariant subspace continuation, e.g., \cite{bindel14}; the goal is to 
track the small eigenvalues in an efficient way, and to check the performance 
of other test functions as an alternative to the determinant of the 
(extended) Jacobian or the number of eigenvalues with negative 
real--parts used so far (see \cite[\S2.1,\S3.1.6]{p2p}). 
\item There are function handles {\tt p.fuha.lss} and {\tt p.fuha.blss} 
for solving the linear systems and the extended linear systems that occur 
in Newton loops and, e.g., for calculating new tangent predictors. 
However, except for some Sherman--Morrison formulas for the case 
of an Allen--Cahn equation with global coupling \cite[\S3.5]{p2p}, 
we always use the standard solvers {\tt lss.m} and {\tt blss.m}, 
which simply call \mlab's $\backslash$ operator. So far, this turned out 
superior to iterative solvers, but this seems to change for 
very large systems, in particular in 3D, see \ref{3D} below, 
and in summary some iterative solvers and customized solvers for 
bordered systems will be fitted into \pdepb\ as well. 
\item A \mlab\ environment online help for \pdepb\ should be 
coming soon. 
\item\label{3D} 
The basic functionality of \pdepb\ has already been ported to 3D, based 
on the free \mlab\ FEM package OOPDE, \cite{uwe}. In particular, this gives 
identical user interfaces in 2D and 3D. In the long term, the 2D, 3D 
(and also 1D) versions 
shall merge to a single package. This is partly similar to the philosophy 
of COCO \cite{coco}. 
\ecen
As \pdep\ is and will remain an ``open project'', comments and help 
on any of the above points will be very welcome. Please send questions, 
remarks or requests to {\tt pde2path@uni-oldenburg.de} or to any of 
the authors.


\renewcommand{\refname}{References}
\renewcommand{\arraystretch}{1.05}\renewcommand{\baselinestretch}{1}
\small
\bibliographystyle{plain}\bibliography{p2p2}

\newpage\appendix
\section{Tables of \pdepb\ functions, controls, switches and fields}
\label{app1}
In this appendix, intended as a reference card, 
we give overviews of the main \pdepb\ functions 
(see the files in {\tt p2plib} for more comments), and of the  
basic \pdepb\ structure {\tt p} and 
the contents of its fields. 
\LTcapwidth=\textwidth
\bce{\small 
\begin{longtable}{|p{15mm}|p{142mm}|} 
\caption{Main fields in the structure {\tt p} describing a p2p2 proplem; see {\tt stanparam.m} in 
{\tt p2plib} for detailed information on the contents of these 
fields and the standard settings. The destinction between {\tt nc} and 
{\tt sw} is somewhat fuzzy, as both contain variables to control 
the behaviour of the numerics: the rule is that {\tt nc} contains 
numerical constants, real or integer, while the switches in {\tt sw} 
only take a finite number of values like 0,1,2,3. Finally, 
{\tt u,np,nu,tau} and {\tt branch} are {\em not} grouped into a substructure 
as, in our experience, these are the variables most often accessed directly 
by the user.\label{tab3}}
\endhead\endfoot\endlastfoot
\hline
field&purpose\\
\hline
fuha&struct of {\bf fu}nction {\bf ha}ndles; 
in particular the function handles 
p.fuha.G, p.fuha.Gjac, p.fuha.bc, p.fuha.bcjac defining 
\reff{gform} and Jacobians, and others such as 
p.fuha.outfu, p.fuha.savefu, ...\\
nc, sw&{\bf n}umerical {\bf c}ontrols such as p.nc.tol, p.nc.nq, \ldots, 
and {\bf sw}itches such as p.sw.bifcheck,\ldots \\[1mm]
u,np,nu& the solution u (including all parameters/auxiliary variables 
in u(p.nu+1:end)), the number of nodes p.np 
in the mesh, and the number of nodal values p.nu of PDE--variables\\
tau,branch&tangent tau(1:p.nu+p.nc.nq+1), and the branch, filled 
via bradat.m and p.fuha.outfu.\\
sol&other values/fields calculated at runtime, e.g.: ds (stepsize), res (residual), \ldots \\
usrlam&vector of user set target 
values for the primary parameter, default usrlam=[];\\ 
eqn,mesh&the tensors $c,a,b$ for the semilinear FEM setup, and 
the geometry data and mesh. \\
plot, file&switches (and, e.g., figure numbers and directory name) 
for plotting and file output\\
time, pm&timing information, and pmcont switches\\
fsol&switches for the interface to {\tt fsolve}, see Remark 2.\\
mat&problem matrices, e.g., mass/stiffness matrices $M$, $K$ for the 
the semilinear FEM setting, and {\tt drop} and {\tt fill} for periodic BC; 
by default, mat is {\em not} saved to disk, see also 
Remark 4.\\
\hline
\end{longtable}
}\ece

{\small 
\bce 
\begin{longtable}{|p{0.26\textwidth}|p{0.68\textwidth}|} 
\caption{Main \pdepb\ functions for user calls; some of these take auxiliary 
parameters, and in general the behaviour is controlled by the settings in 
{\tt p.nc} and  {\tt p.sw}; \ldots indicates additional arguments. 
See the m-files and the 
demo-directories for details.
}\label{mtab}\endfirsthead
\endhead\endfoot\endlastfoot
\hline
function&purpose,remarks\\
\hline
p=stanparam(p)&sets many parameters to ``standard'' values; typically called 
during initialization; also serves as documentation of the meaning 
of parameters\\
p=cont(p), p=pmcont(p)& continuation of problem p, and 
parallel multi-predictor version \\
p=swibra(dir,bptnr,varargin)&branch--switching at point dir/bptnr, varargin 
for new dir and ds\\
plotbra(p,var)&plot branch in p, see also plotbraf.m for plotting from file; 
see also p.plot for settings for plotting\\
plotsol(p,wnr,cmp,style)&plot solution, see also plotsolu, plotsolf, and 
plotEvec\\
p=loadp(dir,pname,varargin)&load p-data at the point pname from directory dir; varargin for new dir\\
p=swipar(p,var)&switch parametrization, see also swiparf\\
p=setpar(p,par)&set parameter values, see also  par=getpar(p,varargin), 
p=setlam(p,lam), and getlam(p);\\
geo=rec(lx,ly)&encode rectangular domain in \ptool\ syntax\\
bc=gnbc(neq,vararg)&generate \ptool--style boundary conditions, see also the 
convenience functions [geo,bc]=recnbc*(lx,ly) and 
[geo,bc]=recdbc*(lx,ly), *=1,2\\
p=findbif(p,varargin)&bifurcation detection via change of stability index; 
alternative to bifurcation detection in cont or pmcont; 
can be run with larger ds, as even number of eigenvalues crossing the imaginary axis is no 
problem\\ 
p=spcontini(dir,name,npar)&initialization for "spectral continuation", 
e.g.~fold continuation\\
p=spcontexit(dir,name)&exit spectral continuation\\
p=rec2per(p)&transform to periodic BC by setting p.mat.drop, p.mat.fill;\\
{}[u,\ldots]=nloop(p,u)&Newton--loop for $(G(u),q(u))=0$\\
{}[u,\ldots]=nloopext(p,u)&Newton--loop for the 
extended system $(G(u),q(u),p(u))=0$\\
p=meshref(p,varargin)&adaptively refine mesh\\
p=meshadac(p)&project onto background mesh p.bmesh, then adaptively refine\\
p=setfemops(p)&set the FEM operators like $M,K$ for the semilinear p.sw.sfem=1 setting\\
p=setfn(p,name)&set output directory to name (or p, if name omitted) \\
err=errcheck(p)&calculate error-estimate\\
screenlayout(p)&position figures for solution-plot, branch-plot and information\\
{}[Gua, Gun]=jaccheck(p)&compare Jacobian p.fuha.Gjac (resp.~p.fuha.sGjac) with 
finite differences\\
p=tint(p,dt,nt,pmod)&time integration of $\pa_t u=-G(u)$; see also tintx for a version with more 
input and output arguments, and saving of selected time-steps.\\%
{}p=tints(p,dt,nt,pmod,nffu)&time integration based on the semilinear
p.sw.sfem=1 setting. If applicable, much faster than tint; again, see also 
tintxs\\
p=loadp2(dir,name,name0)&load u-data from name in directory dir, 
other p-data from name0\\
\hline
\end{longtable}

\begin{longtable}{|p{0.26\textwidth}|p{0.68\textwidth}|}
\caption{Description of functions in {\tt p.fuha}; In the first block, 
only {G, bc, bcjac} are needed if {p.sw.sfem=0}, Gjac (or sGjac) 
only if {p.sw.jac$>$0}. The defaults in the second 
block are set by p=stanparam(p). Third block only 
needed/recommended if p.nc.nq$>0$, 
or for spectral continuation, respectively. 
\label{fuhatab}}
\endfirsthead\endhead\endfoot\endlastfoot
\hline
function&purpose, remarks\\
\hline
{}[c,a,f,b]=G(p,u)&compute coeffcients $c, a, b$ and $f$ in $G$ 
in the full (sfem=0) syntax \\
{}[cj,aj,bj]=Gjac(p,u)&coefficients for calculating $G_u$ in the (sfem=0) 
syntax\\
r=sG(p,u), Gu=sGjac(p,u) & residual $G(u)$ and jacobian $G_u(u)$ 
in the sfem=1 setting 
using the preassembled matrices p.mat.M, p.mat.K, p.mat.Kadv\\
bc=bc(p,u), bcj=bcjac(p,u) &boundary conditions, and their jacobian\\\hline
[p,cstop]=ufu(p,brdat,ds)&user function called 
after each cont.~step, for instance 
to check $\lam_{\min}<\lam<\lam_{\max}$, and to give printout; 
cont.~stops if ufu returns cstop$>$0; default=stanufu, which also checks 
if $\lam$ has passed a value in {\tt p.usrlam}.\\
headfu(p)&function called at start of cont, e.g.~for printout; 
default stanheadfu\\
out=outfu(p,u)&function to generate branch data additional to bradat.m; 
default stanbra\\
savefu(p,varargin)&function to save solution data, default stansavefu; 
see also p.file for settings for saving\\
p=postmmod(p)&function called after mesh-modification; default stanpostmeshmod\\
x=lss(A,u,p)&linear system solver for $Ax=u$, $A=D(G,q)$; default lss with 
$x=A\backslash u$\\
x=blss(A,u,p)&linear system solver for $Ax=u$, $A=D(G,q,p)$ (extended 
or bordered linear system in arclength cont.); default blss with 
$x=A\backslash u$\\ \hline
q=qf(p,u), qu=qjac(p,u)&additional equation(s) $q(u){=}0$, and Jac.~function, 
see, e.g., demo {\tt fCH}\\
Guuphi=spjac(p,u)&$\pa_u(\pa_u G\phi)$ for fold--or branchpoint continuation, 
see, e.g., demo {\tt acfold}\\
\hline
\end{longtable}
\pagebreak

\begin{longtable}{|p{0.21\textwidth}|p{0.75\textwidth}|}
\caption{Description of main numerical controls in {\tt p.nc}.\label{nctab}}
\endfirsthead\endhead\endfoot\endlastfoot
\hline
name and default (where applicable)&purpose, remarks\\
\hline
neq, nq&number $N$ of equations in $G(u)$, see \reff{gform}; 
number of additional equations \reff{qieq}\\
tol=1e-10, imax=10&desired residual; max iterations in Newton loops\\
del=1e-8&stepsize for numerical differentiation\\
ilam&indices of active parameters; ilam(1) is the primary parameter\\ 
lammin,lammax=$\mp 1e6$&bounds for primary parameter during continuation, 
also added to p.usrlam\\
dsmin, dsmax&min and max arclength stepsize, current stepsize in p.sol.ds\\
dsinciter=imax/2&increase ds by factor dsincfac=2 if iter $<$ dsinciter \\
dlammax=1&max stepsize in primary parameter\\
lamdtol=0.5&control to switch between arclength and natural parametrization if 
p.sw.para=1;\\
dsminbis=1e-9&min arclength in bisection for bifurcation localization\\
bisecmax=10&max \# of bisections in bifurcation localization\\
nsteps=10&\# of continuation steps (multiple steps for pmcont)\\
ntot=10000&total maximal \# of continuation steps\\
neig=50&\# of eigenvalues closest to 0 calculated for 
stability (and bif.~in findbif)\\
errbound=0&used as indicator for mesh refinement if $>0$\\
amod=0&mesh-adaption each amod-th step, none if amod=0\\
ngen=3&number of refinement steps under mesh-refinement\\
bddistx=bddisty=0.1&for periodic BC: do not refine at distance$<$ 
bddistx/y from respective boundary\\
\hline
\end{longtable}

\begin{longtable}{| p{0.16\textwidth}|p{0.78\textwidth}|}
\caption{Description of switches in {\tt p.sw}.\label{swtab}}
\endfirsthead\endhead\endfoot\endlastfoot
\hline
name and default&purpose, remarks\\
\hline
bifcheck=1&0/1 for bif.detection off/on\\
spcalc=1&0/1 for calc.~eigenvalues nearest to 0 off/on\\
foldcheck=0&0/1 for fold detection off/on\\
jac=1&0/1 for numerical/analytical (via p.fuha.(s)jac) jacobians for $G$\\
qjac=1&0/1 for numerical/analytical (via p.fuha.qjac) jacobians for $q$\\
spjac=1&0/1 for numerical/analytical (via p.fuha.spjac) jacobian for 
spectral point cont.\\
sfem=0&0/1 for full/semilinear FEM setting\\
newt=0&0/1 for full/chord Newton method\\
bifloc=2&0 for tangent, 1 for secant, 2 for quadratic predictor 
in bif.localization\\
bcper=0&0 for BC via p.fuha.bc, 
1 for  top=bottom, 2 for left=right, 3 for torus\\
spcont=0&0 for normal cont., 1 for bif.~point cont., 2 for fold cont.\\
para=1&0: natural parametr.; 2: arclength; 1: automatic switching
via $\dot\lam<>$p.nc.lamd\\
norm='inf'&or use any number$\ge 1$\\
inter=1,verb=1& interaction and verbosity switches $\in\{0=\text{little},
1=\text{some},2=\text{much}\}$\\
bprint=[]&indices of user-branch data for printout\\
\hline
\end{longtable}

\newpage

\begin{longtable}{| p{0.1\textwidth}|p{0.34\textwidth}| 
p{0.1\textwidth}|p{0.34\textwidth}|}
\caption{Summary of additional data in {\tt p.sol} calculated 
at runtime. Note that the actual solution is stored directly 
in {\tt p.u}, and similarly for the tangent {\tt p.tau}, the 
branch data {\tt p.branch} and the frequently needed data 
{\tt p.np} (number of points in mesh) and {\tt p.nu} (number 
of PDE variables). \label{soltab}}
\endfirsthead\endhead\endfoot\endlastfoot
\hline
name&meaning&name&meaning\\
\hline
deta&sign of det($A$)&muv&vector of eigenvalues of $G_u$\\
err&error estimate&lamd&$\dot\lam$\\
meth&used method (nat or arc)&restart&1 to restart continuation\\
iter&\# of iterations in last Newton loop&xi,xiq&$\xi, \xi_q$ from \reff{xiqeq}\\
ineg&\# of negative eigenvalues&ds&current stepsize\\
\hline
\end{longtable}

\begin{longtable}{| p{0.16\textwidth}|p{0.78\textwidth}|}
\caption{Summary of {\tt p.file}. \label{filetab}}
\endfirsthead\endhead\endfoot\endlastfoot
\hline
name&meaning\\
\hline
count, b(f)count&counters for regular/bif./fold points; 
filenames for regular, bif., fold points 
automatically composed as dir/pt{\tt count}.mat,  
dir/bpt{\tt bcount}.mat and  dir/fpt{\tt fcount}.mat\\
dir, pnamesw=0&directory for saving; 
if pnamesw=1, then set to 'name of p';\\ 
dirchecksw=0& if dirchecksw=1, then warnings given if files 
might be overwritten
\\\hline
\end{longtable}

\begin{longtable}{| p{0.15\textwidth}|p{0.34\textwidth}|p{0.15\textwidth}|p{0.24\textwidth}|}
\caption{Summary of {\tt p.plot}. \label{plottab}}
\endfirsthead\endhead\endfoot\endlastfoot
\hline
name \& default&meaning&name \& default &meaning\\\hline
pfig=1, brfig=2&fig.\,nr.\,for sol./branch plot at runtime& 
ifig=6, spfig=4&info(mesh)/spectrum plot\\
brafig=3&fig.~nr.~for plotbra ( a posteriori) &labelsw=0&axis labeling\\
fs=16&fontsize&lpos=[0 0 10]&light position\\
cm='hot'&colormap&axis='tight'&axis type\\\hline
pstyle=2&\multicolumn{3}{p{0.75\textwidth}|}{plotstyle=0,1,2,3; or customize plotsol}\\
pcmp=1, bpcmp=0&\multicolumn{3}{p{0.75\textwidth}|}
{component\# for sol.~plot and branch plot (relative 
 to data in outfu; last component in bradat=$\|u\|_2$ plotted 
if bpcmp=0)}\\
\hline
\end{longtable}

\begin{longtable}{| p{0.23\textwidth}|p{0.7\textwidth}|}
\caption{Summary of {\tt p.pm} and {\tt p.fsol}. \label{pmtab}}
\endfirsthead\endhead\endfoot\endlastfoot
\hline
name and default&meaning\\\hline
pm: mst=10, imax=1, resfac=0.2 & \# of parallel predictors, \# of iterations in each Newton loop (adapted), factor for desired residual improvement; see 
\cite[\S4.3]{p2p}\\
\hline
fsol: fsol=0, tol=1e-16, imax=5, meth, disp, opt& turn on(1)/off(0) fsol; tol and 
imax for fsol, and {\tt fsolve} options. Note: {\tt fsolve} tolerance 
applies to $\|G(u)\|_2^2$. 
\\
\hline
\end{longtable}

\ece
}

\end{document}